\documentclass[9pt,shortpaper,twoside,web]{ieeecolor}
\usepackage{generic}
\usepackage{cite}
\usepackage{amsmath,amssymb,amsfonts}
\usepackage{algorithmic}
\usepackage{graphicx}

\usepackage{braket}

\usepackage{amsthm}
\usepackage{thmtools}
\usepackage{mathtools}
\usepackage{bigints} 	
\usepackage{steinmetz}  
\usepackage{graphics}  
\usepackage{graphicx}
\let\labelindent\relax
\usepackage{enumitem}
\usepackage{booktabs} 
\usepackage{setspace}
\usepackage{upgreek}
\usepackage{fancyhdr}
\usepackage{float}
\usepackage{tcolorbox}
\usepackage{comment}
\usepackage{empheq, blkarray}
\usepackage{subfig}
\usepackage{subfloat}
\usepackage{array}
\newcolumntype{L}{>{\centering\arraybackslash}m{5cm}}
\makeatletter
\newcommand{\removelatexerror}{\let\@latex@error\@gobble}
\makeatother
\usepackage{verbatim}
\tcbuselibrary{breakable}
\newcommand\numberthis{\addtocounter{equation}{1}\tag{\theequation}}

\theoremstyle{definition}
\newtheorem{definition}{Definition}
\newtheorem{assumption}{Assumption}
\theoremstyle{plain}
\declaretheorem[name=Theorem]{thm}
\declaretheorem[name=Lemma]{lem}
\declaretheorem[name=Proposition]{prop}

\theoremstyle{remark}
\newtheorem{remark}{Remark}
\newtheorem{example}{Example}

\DeclareMathOperator{\dom}{dom}

\DeclareMathOperator{\diag}{diag}

\DeclareMathOperator{\interior}{int}
\DeclareMathOperator{\rinterior}{rint}
\DeclareMathOperator{\rboundary}{rbd}
\DeclareMathOperator{\closure}{cl}
\DeclareMathOperator{\aff}{aff}
\DeclareMathOperator{\range}{ran}

\usepackage{caption}
\renewcommand\thetable{\arabic{table}}

\begin{document}
\title{Continuous-time Discounted Mirror-Descent  Dynamics in Monotone Concave Games}
\author{Bolin Gao and Lacra Pavel%
	\thanks{B. Gao and L. Pavel are with the Department of Electrical and Computer Engineering, University of Toronto, Toronto, ON, M5S 3G4, Canada. Emails:
	{\tt\small bolin.gao@mail.utoronto.ca, pavel@control.utoronto.ca}}%
}

\maketitle

\IEEEpeerreviewmaketitle
\begin{abstract}
	We consider concave continuous-kernel games 
	characterized by monotonicity properties and propose discounted mirror descent-type dynamics. We introduce two classes of dynamics whereby the associated mirror map is constructed based on a strongly convex or a Legendre regularizer. Depending on the properties of the regularizer, 
	we show that these new dynamics can converge asymptotically in concave games with monotone (negative) pseudo-gradient. Furthermore, we show that when the regularizer enjoys strong convexity, the resulting dynamics can converge even in games with hypo-monotone (negative) pseudo-gradient, which corresponds to a shortage of monotonicity. 
\end{abstract}

\section{Introduction}

One of the earliest works on solving continuous-kernel concave games is the work of  Rosen \cite{Rosen65}. The continuous-time gradient type dynamics  was shown to converge to the Nash equilibrium in games that satisfy a so-called diagonally strictly concave condition, roughly equivalent to the pseudo-gradient being a strictly monotone operator. Recently, research on solving monotone games has seen a surge. Both continuous-time dynamics  and discrete-time algorithms have been developed, mostly for games with strictly (strongly) monotone pseudo-gradient. 
For \emph{(non-strictly) monotone only} games, no continuous-time dynamics exist. Discrete-time algorithms  have been proposed, based on either proximal regularization, \cite{Facchinei_I},  inexact proximal best-response,  \cite{Scutari} or Tikhonov type regularization,  \cite{Shanbhag12}, 
and recently extended to generalized Nash equilibrium, e.g. 
 \cite{Peng18},  \cite{GrammaticoDR_CDC2018}. All these works 
 are done in a \emph{discrete-time} setting and the dynamics evolve in the \emph{primal space} of decision variables (and possibly multipliers). With the exception of  \cite{Shanbhag12},  these algorithms 
 are applicable only in games with ``cheap" (inexpensive) proximal/resolvent evaluation, \cite{Facchinei_I}. 
  
  
  In this note we propose a family of \emph{continuous-time discounted mirror descent} dynamics, whereby the dynamics evolves in the space of dual (pseudo-gradient) variables. The mapping from the dual back to the primal space of decision variables is done via a mirror map, constructed based on  two general classes of regularizers. Depending on the properties of the regularizer, we show that these dynamics can converge asymptotically in merely monotone, and even hypo-monotone, concave games.  To the best of our knowledge, these are the first such dynamics in the literature. Our novel contributions consist in relating the convergence of the dynamics to the properties of the convex conjugate of the regularizer.

\emph{Literature review: }  
 Mirror descent algorithms 
  have found numerous applications in recent years, e.g. in distributed optimization \cite{Doan19}, online learning \cite{McMahan11}, and variational inequality problems \cite{Juditsky11}. They fall into the class of so-called \textit{primal-dual} algorithms; the name \textit{mirror descent} refers to the 
two iterative steps: a mapping of the primal variable into a dual space (in the sense of convex conjugate), followed by a mapping of the dual variable, or some post-processing of it, back into the primal space via a mirror map. 
The mirror descent algorithm (MDA) introduced by Nemrovskii and Yudin\cite{Nemirovsky83}, was originally proposed as a generalization of projected gradient descent (PGD) for constrained optimization. The authors of \cite{Beck03} have shown that MDA possesses better rate of convergence as compared to the PGD, which makes it especially suitable for large-scale optimization problems. Other types of algorithms can be seen as equivalent to or special cases of  MDA, e.g. dual averaging \cite{Duchi} and follow-the-leader \cite{McMahan11}. 
A continuous-time version of MDA, 
referred to as the mirror descent  (MD) dynamics, \cite{Staudigl18, Krichene},  
  captures many existing continuous-time dynamics as special cases, such as the gradient flow \cite{Staudigl18, Krichene}, saddle-point dynamics \cite{Cherukuri17} and pseudo-gradient dynamics \cite{Flam02}.  


In the context of  multi-agent games, \emph{mirror descent-like} algorithms have been applied to continuous-kernel games \cite{Zhou}, finite-games \cite{Mertikopoulos16, Bo_LP_CDC2018}, and population games. 
The primal space is taken to be the space of decisions/strategies, and the dual space is the space of payoff vectors (in finite games) or pseudo-gradient vectors (in continuous-kernel games). 
Zhou et al.  \cite{Zhou} introduced the concept of variationally stable concave game 
and showed that,  under variational stability, the iterates of an 
 online MDA converge to the set of Nash equilibria, 
whenever the step-size is  slowly vanishing step-size sequence and that the mirror map satisfies a \textit{Fenchel coupling conforming} condition \cite{Zhou}. Since all concave games with strictly monotone pseudo-gradient  are variationally stable concave games, therefore the algorithm converges in all strictly monotone games. 
However, there are games with a (unique) Nash equilibrium that is not necessarily variationally stable, e.g. zero-sum (monotone) games. While finding the Nash equilibrium of strictly monotone games  
 is an important problem, convergence in such games does not necessarily imply convergence in monotone (but not strictly monotone) games. 
 




\emph{Contributions: } 
Motivated by the above, in this work we propose two  classes of continuos-time \emph{discounted} MD dynamics for concave, continuous-kernel games. The discounting is performed on the dual step of the mirror descent, which generates a \textit{weighted-aggregation} effect similar to the dynamics studied for finite-action games in \cite{Bo_LP_CDC2018}. Discounting is known to foster convergence and eliminate cycling in games, as shown in monotone games or zero-sum games \cite{Bo_LP_CDC2018, Mertikopoulos18}. 
 By exploiting properties of the mirror map in the two classes as well as the discounting effect, we show that these dynamics converge asymptotically to the perturbed equilibria of concave games with monotone (not necessarily strictly monotone) pseudo-gradient. Under certain conditions, they can even converge in concave games with hypo-monotone pseudo-gradient. To the best of our knowledge, these are the first such results. 
Our convergence analysis uses a Lyapunov function given by a Bregman divergence. 
We note that recently \cite{Porco}  
 identified the Bregman divergence as a natural Lyapunov candidate for a variety of systems, elegantly tying with existing results on mirror descent dynamics \cite{Krichene}. 
While the dynamics are in the dual space as in \cite{Bo_LP_CDC2018}, herein we consider continuous-kernel games rather than finite-action games. Furthermore, compared to \cite{Bo_LP_CDC2018} we set up a general framework in terms of two classes of regularizers, matched to the geometry of the action set.  For either strongly convex or Legendre regularizers, we provide convergence guarantees in monotone (hypo-monotone) games   
and present several example discounted MD dynamics. In fact, one such example recovers the dynamics in  \cite{Bo_LP_CDC2018} if the action set is specialized to a simplex  geometry and the regularizer taken as a particular entropy example.  Another example dynamics can be seen as the continuous-time dual counterpart to the  discrete-time Tikhonov (primal) regularization, \cite{Shanbhag12}. Compared to the undiscounted MD \cite{Zhou}, our discounted MD dynamics can converge in (not strictly) monotone games, and even in hypo-monotone games. A short version will appear in \cite{Bo_LP_CDC2019}, with two 
example dynamics. Here we propose two general classes, present proofs (omitted from  \cite{Bo_LP_CDC2019}), additional example dynamics and numerical results. 

The paper is organized as follows. In Section II, we provide preliminary background. 
Section III presents the problem setup and introduces a general form of the discounted mirror descent (DMD) dynamics. 
In Section IV, we construct two classes of DMD and prove their convergence. In Section V, we construct several examples of DMD from each class. 
We present numerical results in Section VI and conclusions in  Section VII. 

\vspace{-0.2cm}
\section{Background}


\vspace{-0.1cm}
\subsection{Convex Sets, Fenchel Duality and Monotone Operators}

The following is from \cite{Beck17, Facchinei_I, Rockafellar}. 
Given a convex set $\mathcal{C} \subseteq \mathbb{R}^n$, the (relative) interior of the set is denoted as ($\rinterior(\mathcal{C})$) $\interior(\mathcal{C})$.  
$\rinterior(\mathcal{C})$ coincides with $\interior(\mathcal{C})$ whenever $\interior(\mathcal{C})$ is non-empty. The closure of $\mathcal{C}$ is denoted as $\closure(\mathcal{C})$, and the relative boundary of $\mathcal{C}$ is defined as $\rboundary(\mathcal{C}) = \closure(\mathcal{C}) \backslash \rinterior(\mathcal{C})$. The indicator function over $\mathcal{C}$ is denoted by $\delta_{\mathcal{C}}$. The normal cone of $\mathcal{C}$ is defined as $N_\Omega(x) =\{v \in \mathbb{R}^n| v^\top(y-x) \leq 0, \forall y \in \mathcal{C}\}$ and  $\pi_{\mathcal{C}}(x) =  \underset{y \in \mathcal{C}}{\text{argmin}} \|y - x\|_2^2$ is the Euclidean projection of $x$ onto $\mathcal{C}$. 


Let $\mathbb{E} = \mathbb{R}^n$ be endowed with norm $\|\cdot\|$ and inner product $\langle \cdot, \cdot \rangle$. An extended real-valued function is a function $f$ that maps from $\mathbb{E}$ to $[-\infty, \infty]$. The (effective) domain  of $f$ is $\dom(f) = \{x \in \mathbb{E} : f(x) < \infty\}$. 
A function $f: \mathbb{E} \to [-\infty, \infty]$ is proper if it does not attain the value $-\infty$ and there exists at least one $x \in \mathbb{E}$ such that $f(x)  < \infty$; it is closed if its epigraph is closed.
A function $f: \mathbb{E} \to [-\infty, \infty]$ is supercoercive if $\lim_{\|x\| \to \infty} f(x)/\|x\| \to \infty$. Let  $\partial f(x)$ denote a subgradient of $f$ at $x$ and  $\nabla f(x)$ the gradient of $f$ at $x$, if $f$ is differentiable.   Suppose $f$ is a closed convex proper on $\mathbb{E}$ with $\interior(\dom(f)) \neq \varnothing$, then $f$ is essentially smooth if $f$ is differentiable on $\interior(\dom(f))$ and $\lim\limits_{k \!\to \!\infty\!} \|\nabla f(x_k)\| \to \infty$ whenever $\{x_k\}_{k = 1}^\infty$ is a sequence in $\interior(\dom(f))$ converging towards a boundary point. 
$f$ 
 is essentially strictly convex if $f$ is strictly convex on every convex subset of  $\dom(\partial f)$. 
A function $f$ is Legendre if it is both essentially smooth and essentially strictly convex.  Given $f$, the function $f^\star \!: \!\mathbb{E}^\star \! \to \! [-\infty, \infty]$ defined by $f^\star(z) \!=\! \sup_{x \! \in \mathbb{E}} \! \big[x^\top z \!- \!f(x)\big]$, 
is called the conjugate function of $f$, where  $\mathbb{E}^\star\!$ is the dual space of $\mathbb{E}$, endowed with the dual norm $\| \cdot\|_\star $.   $f^\star$ is closed and convex if $f$ is proper. By Fenchel's inequality, for any $x \in \mathbb{E}$, $z \in \mathbb{E}^\star$, $f(x) \!+\! f^\star(z)\! \geq \!z^\top x$ {(with equality if and only if $z \! \in \!\partial f(x)$ for proper and convex $f$, or $x \!\in \! \partial f^\star(z)$ if in addition $f$ is closed \cite[Theorem 4.20]{Beck17})}. 
The Bregman divergence of a proper, closed, convex function $f$, differentiable over $\dom(\partial f)$, is $D_f \!:\! \dom(f) \!\times \!\dom(\partial f) \to \mathbb{R}, D_f(x,y) \!=\! f(x) \!-\! f(y) \!-\! \nabla f(y)^\top(x\!-\!y)$. $F \!: \!\mathcal{C} \! \subseteq \!\mathbb{R}^n \!\to \!\mathbb{R}^n$ is monotone  if $(F(z) \!-\! F(z^\prime))^\top (z-z^\prime) \!\geq\! 0$, $\forall z,z^\prime \! \in \! \mathcal{C}$. 
$F$ is $L$-Lipschitz if $\|F(z) \! -\! F(z^\prime)\| \!\leq  \!L \|z \!-\! z^\prime\|$, 
for some $L \!>\! 0$ and is  
\textit{$\beta$-cocoercive} if $(F(z) \!- \!F(z^\prime)^\top (z-z^\prime ) \! \geq \!\beta\|F(z) \!- \!F(z)^\prime\|^2, \forall z, z^\prime \! \in \!\mathcal{C}$ for some $\beta \! > \!0$. 

\vspace{-0.3cm}
\subsection{$N$-Player Concave Games}
Let $\mathcal{G} = (\mathcal{N}, \{\Omega^p\}_{p \in \mathcal{N}}, \{\mathcal{U}^p\}_{p\in \mathcal{N}})$ be a game, where  $\mathcal{N} = \{1, \ldots, N\}$ is the set of players, $\Omega^p \subseteq \mathbb{R}^{n_p}$ is the set of player $p$'s strategies (actions). We denote the strategy (action) set of player $p$'s opponents as $\Omega^{-p} \subseteq \prod_{q \in \mathcal{N}, q \neq p} \mathbb{R}^{n_q}$. We denote the set of all the players strategies as $\Omega =   \prod_{p \in \mathcal{N}} \Omega^{p} \subseteq \prod_{p \in \mathcal{N}} \mathbb{R}^{n_p} = \mathbb{R}^{n}, n = \sum_{p \in \mathcal{N}} n_p$. We refer to  $\mathcal{U}^p: \Omega \to \mathbb{R}, x \mapsto \mathcal{U}^p(x)$ as player $p$'s real-valued payoff function, where $x  = (x^p)_{p \in \mathcal{N}} \in \Omega$ is the action profile of all players, and $x^p \in \Omega^p$ is the action of player $p$. We also  denote $x$ as $x  = (x^p;x^{-p})$ 
where $x^{-p} \in \Omega^{-p}$ is the action profile of all players except $p$. 
\vspace{-0.23cm}
\begin{assumption} For all $p \in \mathcal{N}$,
	\label{assump:blanket}
	\begin{itemize}
		\item[i.]  $\Omega^p$ is a non-empty, closed, convex, subset of $\mathbb{R}^{n_p}$,
		\item[ii.] $\mathcal{U}^p(x^p;x^{-p})$ is (jointly) continuous in $x = (x^p;x^{-p})$,
		\item[iii.] $\mathcal{U}^p(x^p;x^{-p})$ is concave and continuously differentiable in each $x^p$ for all $x^{-p} \in \Omega^{-p}$. 
	\end{itemize}
\end{assumption} \vspace{-0.23cm}
Under \autoref{assump:blanket}, we refer to $\mathcal{G}$ as a \textit{concave game}. Equivalently, in terms of a cost function  $J^p \!= \!-\mathcal{U}^p$, the game $\mathcal{G}$ is a \textit{convex game}. For the rest of the paper, we use the payoff function throughout. 
Given $x^{-p} \in \Omega^p$, 
each agent $p \in \mathcal{N}$ aims to find the solution of the following optimization problem,  \vspace{-0.23cm}
\begin{equation}
\begin{aligned}
& \underset{x^p}{\text{maximize}}
& &  \mathcal{U}^p(x^p; x^{-p})
& \text{subject to}
& & x^p \in \Omega^p.
\end{aligned}
\end{equation}
A profile 	${{x}}^\star \!=\! ({{x}^p}^\star)_{p \in \mathcal{N}} \!\in \!\Omega^p$  is a Nash equilibrium if, \vspace{-0.2cm}
	\begin{equation}\label{Nash_definition}
	\mathcal{U}^p({x^p}^\star; {x^{-p}}^\star) \geq \mathcal{U}^p(x^p; {x^{-p}}^\star), \forall x^p \in \Omega^p, \forall p \in \mathcal{N}.
	\end{equation}
At a Nash equilibrium, no player can increase his payoff by unilateral deviation. If $\Omega^p$  is bounded, under Assumption \ref{assump:blanket}, existence of a Nash equilibrium is guaranteed (cf., e.g. \cite[Theorem 4.4]{Basar}). 
When $\Omega^p$ is closed but not bounded, 
existence of a Nash equilibrium is guaranteed under the additional assumption that  $-\mathcal{U}^p$ is coercive in $x^p$, that is, 
$	\lim\limits_{\|x^p\|\to\infty} -\mathcal{U}^p(x^p;x^{-p}) = +\infty,$ 
	for all $x^{-p} \in \Omega^{-p}, p \in \mathcal{N}$, (cf.  \cite[Corollary 4.2]{Basar}).
A useful characterization of a Nash equilibrium of a concave game $\mathcal{G}$ is given in terms of the \textit{pseudo-gradient} defined  as $U: \Omega \to \mathbb{R}^n, U(x) \!=\!(U^p(x))_{p \in \mathcal{N}}$, where 
$U^p(x) =\nabla_{x^p} \mathcal{U}^p(x^p;x^{-p})$ is the partial-gradient. 
By \cite[Proposition 1.4.2] {Facchinei_I}, $x^\star \in \Omega$ is a Nash equilibrium if and only if,\vspace{-0.2cm}
	\begin{equation}
	(x-x^{\star} )^\top U(x^{\star}) \leq 0, \forall x \in \Omega. 
	\label{eqn:nash_equilibrium_VI}
	\end{equation}	
Equivalently $x^\star$ is a solution 
of the variational inequality 
$\text{VI}(\Omega, -U)$,  \cite{Facchinei_I}, or, using the definition of the normal cone, 
\vspace{-0.23cm}
	\begin{equation}\label{NE_equiv}
	U(x^\star) \in N_\Omega(x^\star). 
	\end{equation}

Standard assumptions on the pseudo-gradient are as follows.\vspace{-0.23cm}
\begin{assumption}
	\label{assump:pseudo_gradient}
	$-U(x) =-(\nabla_{x^p} \mathcal{U}^p(x^p, x^{-p}))_{p \in \mathcal{N}}$ is \vspace{-0.1cm}
	\begin{itemize}
		\item [(i)] monotone, 
$		-(U(x) - U(x^\prime))^\top(x-x^\prime) \geq 0, \forall x, x^\prime \in \Omega. 
$
	\item[(ii)] strictly monotone, 
$		-(U(x) \!-\! U(x^\prime))^{\!\top}\!(x\!-\!x^\prime) \!>\! 0, \forall x \! \neq \!x^\prime  \in \Omega.
$
	\item[(iii)] $\eta$-strongly monotone, 
$		\!-(U(x) \!-\! U(x^\prime\!))\!^\top\!(x\!-\!x^\prime) \!\geq \! \eta\|x \!-\! x^\prime\|^2_2$, $\! \forall x,x^\prime  \!\in \!\Omega$, for some $\eta \!>\! 0$.
	\item[(iv)] $\mu$-hypo monotone, 
$		\!-(U(x) \!-\! U(x^\prime\!))\!^\top\!(x\!-\!x^\prime) \!\geq \! -\mu\|x \!-\! x^\prime\|^2_2$, $\! \forall x,x^\prime  \!\in \!\Omega$, for some $\mu \!>\! 0$.
	\end{itemize} 
\end{assumption}\vspace{-0.25cm}
We refer to $\mathcal{G}$ as a monotone game if it satisfies \autoref{assump:pseudo_gradient}(i). 

\vspace{-0.25cm}
\section{Problem Setup}

We consider a set of players who are repeatedly interacting in a concave game $\mathcal{G}$. Assume that the game repeats with an infinitesimal time-step between each stage, hence we model it as a continuous-time process as in \cite{MertikCDC2017}, \cite{Flam02}. 
Each player maps his own partial-gradient $u^p \!=\! U^p(x) \!\in \! \mathbb{R}^{n_p}$ into an \textit{auxiliary variable} $z^p \! \in \! \mathbb{R}^{n_p}$ via a dynamical system $\dot{z}^p \!=\! F(z^p,u^p)$ and selects the next action $x^p \!\in\! \Omega^p$ via a so-called \textit{mirror map} $C^p$.  
The entire learning process for each player can be written as a continuous-time dynamical system,\vspace{-0.25cm}   
\begin{equation}
\label{eqn:general_learning_process}
\begin{cases}
\dot z^p &= F(z^p,u^p),\\
x^p &= C^p(z^p),
\end{cases}
\end{equation}
where $u^p = U^p(x) = \nabla_{x^p} \mathcal{U}^p(x^p;x^{-p})$. We assume that the mirror map $C^p: \mathbb{R}^{n_p}  \to \Omega^p$  is given by,\vspace{-0.23cm}
\begin{equation} 
	\label{eqn:mirror_map_argmax_char} 
	C^p(z^p) = \underset{y^p \in \Omega^p}{\text{argmax}}\left[ {y^p}^\top z^p - \epsilon\vartheta^p(y^p)\right], 
\end{equation} 
where $\vartheta^p: \mathbb{R}^{n_p} \to \mathbb{R}^{n_p}\cup\{\infty\}$ is assumed to be a closed, proper and (at-least) essentially strictly convex function, where $\dom(\vartheta^p) = \Omega^p$ is assumed be a non-empty, closed and convex set. The function $\vartheta^p$ is often referred to as a \emph{regularizer} in optimization, learning and game contexts. Different forms of mirror map can be derived depending on the regularizer. 
Finally, since the pseudo-gradient is not assumed to be bounded, $\vartheta^p$ should be chosen so that the dual space is unconstrained.

The most important family of algorithms that follows the model of the learning dynamics  \eqref{eqn:general_learning_process} is that of mirror descent (MD) dynamics,\vspace{-0.25cm}  
\begin{equation}
\label{eqn:mirror_descent} 
\begin{cases}
\dot z^p &= \gamma u^p,\\
x^p &= C^p(z^p),
\end{cases}
\end{equation}
where 
$\gamma \!>\! 0$ is a rate parameter. 
This 
can be interpreted as each player performing an aggregation of its own partial-gradient, $z^p(t) = z^p(0) +\gamma \int\limits_0^t u^p(\tau) \mathrm{d}\tau$, and mapping it to an action via the mirror map $C^p$. 
The discrete-time analog of \eqref{eqn:mirror_descent}, \vspace{-0.25cm}
\begin{equation}
\begin{cases}
z^p_{k+1} &= z^p_k + \gamma t_k u_k^p,\\
x^p_{k+1} &= C^p(z^p_{k+1}),
\end{cases}
\end{equation} with $\!t_k \!>\! 0$ the step-size,  is the online mirror descent  studied in \cite{Zhou} in  a similar concave game setup. In finite games, this algorithm is referred to as Follow-the-Regularized-Leader (FTRL)\cite{Mertikopoulos18}. 

\begin{remark} 
As an example, let $\vartheta^p(x^p) \!=\! \dfrac{1}{2\epsilon}\|x^p\|^2_2$,  $\dom(\vartheta^p)\! =\! \mathbb{R}^{n_p}$, so, cf. \eqref{eqn:mirror_map_argmax_char}, 
$C^p(z^p) \!=\!  \underset{y^p \in \Omega^p}{\text{argmax}} \left[ {y^p}^\top z^p \!-\!\dfrac{1}{2}\|x^p\|^2_2\right] \!=\! z^p$. The  \emph{dual (MD) dynamics} \eqref{eqn:mirror_descent} is,  \vspace{-0.25cm}
\begin{equation}
 \label{eqn:PSGD-Dual}
\dot z^p 
= \gamma u^p, \quad 
x^p 
= z^p,
\end{equation} which is in turn equivalent to the well-known \emph{primal dynamics},  \vspace{-0.2cm} \begin{equation} \label{eqn:PSGD} \dot x^p = \gamma u^p =  \gamma \nabla_{x^p} \mathcal{U}^p(x^p;x^{-p}), \end{equation} 
or the \emph{pseudo-gradient dynamics} (PSGD),  known to converge to the NE when $-U(x)$ is  strictly/strongly monotone (e.g. Lemma 2, \cite{Dian_Pavel_TAC2019}). 
\end{remark}
\vspace{-0.1cm}
In this paper we propose 
a related variant of the  MD dynamics  \eqref{eqn:mirror_descent}, called the \emph{discounted mirror descent dynamics} DMD, given by,\vspace{-0.25cm}
\begin{equation}
\label{eqn:DMD} 
\begin{cases}
\dot z^p &= \gamma (-z^p + u^p),\\
x^p &= C^p(z^p),\qquad  \quad \forall p \in \mathcal{N}
\end{cases}
\end{equation}
where $u^p \!=\! U^p(x) \!= \!\nabla_{x^p} \mathcal{U}^p(x^p;x^{-p})$, and  $\gamma \!>\! 0$.  
Unlike the undiscounted MD \eqref{eqn:mirror_descent}, in \eqref{eqn:DMD} each player performs an \emph{exponentially discounted} aggregation. 
%
The DMD dynamics of all players can be written  in stacked notation as,\vspace{-0.25cm}
\begin{align}
\dot z & = \gamma(-z+ u), \quad 
x = C(z)
\label{eqn:DMD_overall}
\end{align}
with $x\! \in \! \Omega, z \!\in \!\mathbb{R}^n\!, \!u \!=\!  U(x) \!=\!(U^p(x))_{p \in \mathcal{N}}$, $\!C(z) \!=\! (C^p(z^p))_{p \in \mathcal{N}}$. 

Our focus in this paper is to construct classes of DMD dynamics \eqref{eqn:DMD}  for different types of the regularizer $\vartheta^p$,  \eqref{eqn:mirror_map_argmax_char}. We investigate the convergence of these classes of dynamics in monotone  (not necessarily strictly monotone) games, 
based on the properties of the associated mirror map $C^p$, \eqref{eqn:mirror_map_argmax_char}. We then construct several examples of DMD dynamics from each class. 
\vspace{-0.3cm}
\section{A General Framework for Designing Discounted Mirror Descent Dynamics}

In this section, we consider two general classes of regularizers and study properties of the associated mirror maps (proofs are given in the Appendix). Based on these, we investigate the convergence of DMD  \eqref{eqn:DMD}, under different assumptions on the game's pseudo-gradient. 

\vspace{-0.3cm}
\subsection{Properties of Induced Mirror Maps}

We consider convex regularizers that can be classified as either steep or non-steep according to the following definition.
\vspace{-0.2cm}
\begin{definition}
	A closed, proper, convex regularizer $\vartheta^p\!\!:\! \mathbb{R}^{n_p}\! \!\to\!  \mathbb{R} \cup\! \{\infty\!\}$  is said to be \textit{steep} (or \textit{relatively essentially smooth}) if,
	\begin{itemize}
		\item[(i)]  $\dom(\vartheta^p)$ is non-empty and convex, 
		\item[(ii)] $\vartheta^p$ is differentiable on $\rinterior(\dom(\vartheta^p))$, 
		\item[(iii)] $	\lim\limits_{k \to \infty} \|\nabla \vartheta^p(x^p_k)\| \! =\!  + \infty$, 
		whenever $\{x^p_k\}_{k = 1}^\infty$ is a sequence in $\rinterior(\dom(\vartheta^p))$ converging to a point in $\rboundary(\dom(\vartheta^p))$. 
	\end{itemize} 
	$\vartheta^p$  is \textit{non-steep} if $\lim\limits_{k \to \infty} \|\nabla \vartheta^p(x^p_k)\|$ is bounded, for any sequence $\{x^p_k\}_{k = 1}^\infty$ in $\rinterior(\dom(\vartheta^p))$ converging to a point in $ \rboundary(\dom(\vartheta^p))$.  
\end{definition}

\vspace{-0.3cm}\begin{remark}A non-empty, convex domain $\dom(\vartheta^p)$ ensures the non-emptiness of its relative interior \cite[Theorem 6.2, p. 45]{Rockafellar}. 
\end{remark}

\vspace{-0.3cm}
\begin{prop}
	\label{prop:steep_implies_differentiability_in_interior}
	Let $\vartheta^p\!:\! \mathbb{R}^{n_p} \!\to \!\mathbb{R}\!\cup\!\{\infty\}$ be a closed, proper, convex. 
		Then, the following hold: 
		(i) \!If $\!\vartheta^p\!\!$ is steep, then $\! \!\rboundary(\dom(\vartheta^p\!)\!) \!\not\subset \! \!\dom(\partial \vartheta^p\!)$ and 
		$\!\dom(\partial \vartheta^p\!) \!\!= \!\rinterior(\dom(\vartheta^p\!))$. 
		
		(ii) \! If $\vartheta^p$ is non-steep, then $\!\rboundary(\dom(\vartheta^p))\! \subset \!\!\dom(\partial\vartheta^p)\!$ and 
		$\!\dom(\partial \vartheta^p)\! \!= \!\dom(\vartheta^p\!)$.  
\end{prop}

\vspace{-0.3cm}
\begin{assumption}
	\label{assump:primal}
	The regularizer $\vartheta^p\!:\! \mathbb{R}^{n_p} \!\!\to \!\mathbb{R}\!\cup\!\{\!\infty \!\}$ is closed, proper, convex, with $\dom(\vartheta^p)$ non-empty, closed and convex. In addition, \vspace{-0.2cm}
	\begin{itemize}
		\item[(i)] $\vartheta^p$ is $\rho$-strongly convex,  or
		\item[(ii)] $\vartheta^p$ is Legendre and $\interior(\dom(\vartheta^p))\! \neq \!\varnothing$.
	\end{itemize}
\end{assumption}
\vspace{-0.2cm}
Note that \autoref{assump:primal}(ii) relaxes strong convexity to essential strict convexity and essential smoothness (steepness).  
In order to take into consideration $\epsilon$ in the regularization, 
	cf. \!\eqref{eqn:mirror_map_argmax_char}, we consider $\!\psi^p \! =\!\epsilon \vartheta^p\!$, which inherits all properties of $\vartheta^p\!$.  
	We then refer to $C^p$ as the mirror map induced by $\psi^p\!$. Next, we derive properties of  $C^p\!$ for the two classes of regularizers cf. \autoref{assump:primal}(i) and \autoref{assump:primal}(ii).  
\vspace{-0.2cm}
\begin{prop}
	\label{prop:primal_sc}
	Let $\psi^p\! \!= \!\epsilon \vartheta^p\!$, $\epsilon \!>\!\! 0$, where  $\vartheta^p$ satisfies \autoref{assump:primal}(i), and let ${\psi^p}^\star$ be the convex conjugate of $\psi^p$. 
	Then, 
	\begin{enumerate}
		\item[(i)] ${\psi^p}^\star\!:\!\mathbb{R}^{n_p} \!  \to \!\mathbb{R}\cup\{\infty\}$ is closed, proper, convex and finite-valued over $\mathbb{R}^{n_p}$, i.e., $\dom({\psi^p}^\star) = \mathbb{R}^{n_p}$.
		\item[(ii)]  ${\psi^p}^\star$ is continuously differentiable on $ \mathbb{R}^{n_p}$ and $\nabla {\psi^p}^\star \!= \!C^p$. 
		\item[(iii)] $C^p$ is  $(\epsilon\rho)^{-1}$-Lipschitz on $ \mathbb{R}^{n_p}$. \
		\item[(iv)]  $C^p$ is $\epsilon\rho$-cocoercive on $ \mathbb{R}^{n_p}$, and in particular, is monotone. 
		\item[(v)]  $C^p\!$ is surjective from $\mathbb{R}^{n_p}\!$ onto $\rinterior(\dom(\psi^p\!))$ whenever $\psi^p\!$ is steep, and  onto $\dom(\psi^p\!)$ whenever $\psi^p\!$ is non-steep.  
		\item[(vi)] $C^p\!$ has $\nabla \psi^p\!$ as a left-inverse  over $\rinterior(\dom(\psi^p\!))$ whenever $\psi^p$ is steep, 
		 and  over $\dom(\psi^p)$ whenever $\psi^p\!$ is non-steep.  
	\end{enumerate}
\end{prop}


\vspace{-0.3cm}
\begin{remark}
	If $\psi^p \!=\! \epsilon \vartheta^p\!$ is differentiable over all $\dom(\psi^p\!)$, following \cite[Theorem 6.2.4(b), p. 264]{Niculescu}, \autoref{prop:primal_sc} strengthens as follows:
(i)		${\psi^p}^\star\!$ is closed, proper, strictly convex and finite-valued over $\mathbb{R}^{n_p}\!$, 
(ii)		$C^p$ is strictly monotone on $\mathbb{R}^{n_p}$, 
(iii) $C^p$ is bijective from $\mathbb{R}^{n_p}$ to $\dom(\psi^p)$, 
(iv) $C^p$ has a full inverse $\nabla \psi^p$ over $\dom(\psi^p)$. 
	For example, $\vartheta^p(x^p) \! =\! \dfrac{1}{2\epsilon} \|x^p\|_2^2, x^p \!\in \! \mathbb{R}^{n_p}$, (PSGD) is such a case. 
\end{remark}

\vspace{-0.2cm}
\begin{prop}
	\label{prop:primal_Legendre}
Let $\psi^p\! \!= \!\epsilon \vartheta^p\!$, $\epsilon \!>\!\! 0$, where  $\vartheta^p$ satisfies \autoref{assump:primal}(ii), and let ${\psi^p}^\star$ be the convex conjugate of $\psi^p$. 
Then, 
	\begin{itemize}
		\item[(i)]  ${\psi^p}^\star\!:\! \mathbb{R}^{n_p} \!\to \!\mathbb{R}\!\cup\!\{\infty\}$ is closed, proper, Legendre and finite-valued over $\mathbb{R}^{n_p}$, i.e., $\dom({\psi^p}^\star) \!=\! \mathbb{R}^{n_p}$.
		\item[(ii)]  	$\nabla \psi^p\! :\!  \interior(\dom(\psi^p))\!  \to\!  \interior(\dom({\psi^p}^\star))$ is a homeomorphism with inverse mapping $(\nabla \psi^p)^{-1} \!=\! \nabla {\psi^p}^\star \!=\! C^p.$ 
		\item[(iii)]  $C^p$ is strictly monotone on $\interior(\dom({\psi^p}^\star))$. 
	\end{itemize}	
\end{prop}
\vspace{-0.2cm}
\noindent \autoref{prop:primal_Legendre} follows from Legendre theorem \cite[Thm 26.5, p.258]{Rockafellar}. 

Next, we provide a fixed-point characterization of the mirror map $C^p$ (\autoref{prop:mirrormap_sc}), which will be used to relate equilibria of   \eqref{eqn:DMD_overall} to Nash equilibria of the game (\autoref{prop:pertNash}).
\vspace{-0.2cm}
\begin{prop}\label{prop:mirrormap_sc}
Let $\psi^p\! \!= \!\epsilon \vartheta^p\!$, $\epsilon \!>\!\! 0$, where  $\vartheta^p$ satisfies \autoref{assump:primal}. 
Then, the  mirror map induced by $\psi^p$, $C^p$, \eqref{eqn:mirror_map_argmax_char},  
  can be written as the fixed point of the Bregman projection, \vspace{-0.2cm}
	\begin{equation}
	\label{eqn:Bregman_projection_sc}
	C^p(z^p) =\underset{y^p \in \Omega^p}{\text{argmin}} \thinspace D_{\psi^p}(y^p, C^p(z^p)), \vspace{-0.2cm}
	\end{equation}
	where $D_{\psi^p}$ is the Bregman divergence of $\psi^p$, \vspace{-0.2cm}
	$$D_{\psi^p} (y^p,q^p) \!=\! \psi^p(y^p) \!-\! \psi^p(q^p) \!-\! \nabla \psi^p(q^p)^\top(y^p\!-\!q^p).$$
\end{prop} 
\vspace{-0.2cm}

We show next that any rest point $\overline x$ of DMD \eqref{eqn:DMD} or \eqref{eqn:DMD_overall}  is the  Nash equilibrium associated with a perturbed payoff. 
Any equilibrium point of the closed-loop system   \eqref{eqn:DMD_overall} is characterized by, \vspace{-0.26cm}
\begin{equation}
\label{eqn:DMD_rest_points}
\overline{u}  = U(\overline x) = \overline z, \qquad 
\overline x = C(\overline z),
\end{equation}
i.e.,  $\!\overline z \!=\!  \!U \!\circ \!C(\overline z)$, $\overline x \!= \! \!C\circ U(\overline x)$.  From 
\eqref{eqn:mirror_map_argmax_char},  
by Berge's maximum theorem, 
$C$ is compact valued and upper semicontinuous. Since $U$ is jointly continuous, 
$U \circ C\!$ is also compact and upper semicontinuous, and by Kakutani's fixed-point theorem, admits a fixed point.
\vspace{-0.23cm}
\begin{prop}\label{prop:pertNash}
Let $\psi^p\! \!= \!\epsilon \vartheta^p\!$, $\!\epsilon \!>\!\! 0$, where  $\vartheta^p$ satisfies \autoref{assump:primal} and $C^p\!$ the induced mirror map. 
	Any rest point $\!\overline x \!=\! C(\overline z)$ of  DMD \eqref{eqn:DMD} 
	is the Nash equilibrium of the game $\mathcal{G}$ with perturbed payoff, 
	\vspace{-0.2cm}
	\begin{equation} \label{eqn:perturbed_payoff} \widetilde{\mathcal{U}}^p(x^p; x^{-p}) = \mathcal{U}^p(x^p; x^{-p})  - \epsilon\vartheta^p(x^p).\end{equation} As $\epsilon \!\to \! 0$, $\overline x  \!\to \!x^\star$, where $x^\star \!=\! ({x^p}^\star)_{p = 1}^N$ is a Nash equilibrium of $\mathcal{G}$. 
\end{prop}\vspace{-0.23cm}
\begin{proof}
	From the fixed-point characterization of the mirror map  \eqref{eqn:Bregman_projection_sc} (cf. \autoref{prop:mirrormap_sc}), evaluated at $\overline z =\overline u$, one can write $\forall p$,
\vspace{-0.2cm}
	\begin{align*}
	C^p(\overline{u}^p) = \underset{y^p \in \mathbb{R}^{n_p} }{\text{argmin}} \thinspace \left[ \delta_{\Omega^p}(y^p) + D_{\psi^p}(y^p, C^p(\overline{u}^p))  \right], 
	\end{align*} 
	where $\delta_{\Omega^p}(y^p) $ is the indicator function over  $\Omega^p$. By Fermat's condition for unconstrained optimization \cite[Prop 27.1, p. 497]{Bauschke2nd}, $\overline x^p = C^p(\overline{z}^p) =  C^p(\overline{u}^p)$ is a minimizer if and only if, \vspace{-0.2cm}
	\begin{align} &\mathbf{0} \in \partial \delta_{\Omega^p}(\overline x^p) + \nabla_{y^p}  D_{\psi^p}(\overline x^p,   C^p(\overline{u}^p)), 
	\end{align}
	or 
	$\qquad \mathbf{0} \in  N_{\Omega^p}(\overline x^p) +  \nabla \psi^p(\overline x^p) - \nabla \psi^p(C^p(\overline u^p))$, where    $\partial \delta_{\Omega^p}(\overline x^p) = N_{\Omega^p}(\overline x^p)$ \cite{Bauschke2nd} was used. 
	By \autoref{prop:primal_sc}(ii) or \autoref{prop:primal_Legendre}(ii), $C^p(\overline{u}^p) \!=\! \nabla {\psi^p}^\star(\overline u^p)$, and $C^p$ has  $\nabla \psi^p$  as a left-inverse (cf. \autoref{prop:primal_sc}(vi) or \autoref{prop:primal_Legendre}(ii)),   therefore, $\nabla \psi^p(C^p(\overline u^p)) \!=\!\nabla \psi^p(\nabla {\psi^p}^\star(\overline u^p)) \!=\! \overline u^p$. Substituting this and $\overline u^p \!=\!\nabla_{x^p} \mathcal{U}^p(\overline{x}^p; {\overline{x}^{-p}})$ 
yields  for any $\epsilon\!>\!0$, $\forall p$,\vspace{-0.2cm}
	\[
		\nabla_{ x^p} \mathcal{U}^p(\overline x^p, \overline x^{-p})  \in (N_{\Omega^p} + \nabla \psi^p)(\overline x^p) =  (N_{\Omega^p} + \epsilon \nabla \vartheta^p)(\overline x^p),
	\]
	In stacked form, with $\nabla \vartheta(\overline x)\!=\!(\nabla \vartheta^p(\overline x^p))_{p \in \mathcal{N}}$, this is written as\vspace{-0.2cm}
	\begin{equation}\label{eq:pNE}
	U(\overline x) -  \epsilon   \nabla \vartheta(\overline x) \in   N_{\Omega}(\overline x),
	\end{equation}
or  
$	\!\nabla_{x^p} \widetilde{\mathcal{U}}^p(x^p; \!x^{-p}))_{p = 1}^N \! \in  \! N_{\Omega}\!(\overline x)\! $. By \eqref{NE_equiv}, $\overline x$ is a Nash equilibrium for the perturbed payoff $\widetilde{\mathcal{U}}^p\!$. 
As $\!\epsilon \! \to \! 0$, $\!$\eqref{eq:pNE} yields $\!$ \eqref{NE_equiv}, 
	hence $\! \overline x \!\to \! {x^\star}\!$. 
\end{proof}

\vspace{-0.3cm}
\begin{remark} If $-U$ is monotone, then $\!-(U(x) \!-\!  \epsilon \!\nabla \vartheta(x)\!)\!$ is strictly monotone, hence a unique perturbed NE exists for each $\epsilon \!>\!0$. 
\end{remark}

\vspace{-0.3cm}
\subsection{Convergence of DMD under Induced Mirror Maps}

Using key properties given by \autoref{prop:primal_sc} and \autoref{prop:primal_Legendre}, for regularizers satisfying either  \autoref{assump:primal}(i) or \autoref{assump:primal}(ii), in \autoref{thm:convergence_sc} and \autoref{thm:convergence_Legendre} we show convergence of DMD under corresponding induced mirror maps in the two cases, respectively. \vspace{-0.2cm}
\begin{thm}   
	\label{thm:convergence_sc}
	Let $\mathcal{G} \!= \!(\mathcal{N}, (\Omega^p)_{p \in \mathcal{N}}, (\mathcal{U}^p)_{p \! \in \!\mathcal{N}})$ be a concave game with players' dynamics given by DMD \eqref{eqn:DMD}. Assume there are a finite number of isolated fixed-points $\overline z\!$ of $U \circ C\!$, where $C \!=\! (C^p)_{p \in \mathcal{N}}\!$ is the mirror map induced by $\psi^p \!=\! \epsilon \vartheta^p\!$ satisfying \autoref{assump:primal}(i). Then, under either \autoref{assump:pseudo_gradient}(i), (ii), or (iii),  with the additional assumption that $-\mathcal{U}^p$ is coercive in $x^p$ whenever $\Omega^p\!$ is non-compact, for any $\epsilon \!>\! 0$,  the auxiliary variables $z(t) \!=\! (z^p(t))_{p \in \mathcal{N}}\!$ converge to a rest point $\overline z$ while players' actions $x(t) \!=\! (x^p(t))_{p \in \mathcal{N}}\!$ converge to $\overline x$, a perturbed Nash equilibrium of $\mathcal{G}$. Alternatively, under \autoref{assump:pseudo_gradient}(iv), the same conclusions hold for any $\epsilon \!>\! \mu\rho^{-1}$.  
\end{thm}
\vspace{-0.3cm}
\begin{proof}  
Let  $\overline z$ be a rest point of  \eqref{eqn:DMD_overall}, $\overline z \!=\!\overline{u}\!=\! U(C(\overline z))$. Take as Lyapunov function 
 the sum of Bregman divergences of ${\psi^p}^\star$, 
 $V(z) \!=\! 	\sum_{p \in \mathcal{N}} D_{{\psi^p}^\star}(z^p, \overline z^p)$,  \vspace{-0.2cm}
	\begin{align} \! V(z) &\! =\! \sum\limits_{p \in \mathcal{N}}  {\psi^p}^\star(z^p) \!-\! {\psi^p}^\star(\bar z^p) \!-\! \nabla {\psi^p}^\star({\overline z}^p)^\top(z^p\! -\! {\overline z}^p).  \label{eqn:global_convegence_bregman_sc} \end{align} 
	Since ${\psi^p}^\star$ is convex (cf. \autoref{prop:primal_sc}(i)), it follows that $V$ is  positive semidefinite. 
	When $\Omega^p$ is compact, since $U^p_i(x^p; x^{-p})$ is continuous,  $|U^p_i(x^p; x^{-p})| \!\leq \!M, \forall x  \!\in \!\Omega$, for some $M\!>\!0$. 
	Then from \eqref{eqn:DMD},  $\!|z^p_i\!(t)\!| 
	\! \leq  \! e^{-\gamma t}\! |z^p_i\!(0)\!| \!+\! M\!(1\!-\!e^{-\gamma t}\!)$, and $|z^p_i(t)|\! \leq \! \max\{z_i^p\!(0)\!, M\!\}, \!\forall t \!\geq\! 0.$ Hence $\mathcal{D} \!=\! \{z \!\in \!\mathbb{R}^{n}| \|z\|_2  \!\leq \!\sqrt{n}M\}$ is nonempty, 
	 compact, positively invariant set. Alternatively, when $\Omega^p$ is non-compact, for any $\bar x^p \! \!\in\!\! \interior(\dom(\psi^p))$, ${\psi^p}^\star(\cdot) -\! < \!\bar x^p,\cdot\!>$ is  coercive \cite[Prop. 1.3.9(i)]{Hiriart-Urruty},  hence $V$ is coercive and $\mathcal{D}$ can be any of its sublevel sets. 	
	Along any solution of \eqref{eqn:DMD}, $\dot V(z) \! =\! \sum\limits_{p \in \mathcal{N}} \nabla_{z^p} {D_{{\psi^p}^\star}(z^p, \overline z^p)}^\top \dot z^p 
	\!=\! \sum\limits_{p \in \mathcal{N}}  \gamma(\nabla{\psi^p}^\star(z^p) \! -\! \nabla{\psi^p}^\star({\overline z }^p))^\top(u^p \!-\! z^p\!) $. Using \autoref{prop:primal_sc}(ii),   \vspace{-0.22cm}
	\begin{align*}
		\dot V(z) 
	&	= \sum\limits_{p \in \mathcal{N}} \gamma  (C^p(z^p) - C^p({\overline z}^p))^\top(u^p  - \overline{u}^p + \overline z^p\! -\! z^p)  \nonumber  \\
	&= \gamma (x\! -\! \overline x)^\top(u - \overline{u}) \!-\! \gamma( C(z)\!-\! C(\overline z))^\top(z \!-\! \overline z),  \numberthis \label{eqn:strongly_convex_proof_1} 	\end{align*}
	where 
	$x\! = \! C(z)$ and $\overline x \!= \!C(\overline z)$, cf.\eqref{eqn:DMD_rest_points} was used. 
	Since $u \!= \!U(x), \overline u\! =\! U(\overline x)$, under \autoref{assump:pseudo_gradient}(i), \ref{assump:pseudo_gradient}(ii), or \ref{assump:pseudo_gradient}(iii) the first term of $\dot V(z)$ 
	is non-positive, therefore,
	$	\dot V(z) \! \leq \! -\gamma (C(z)-C(\overline z))^\top(z-\overline z)
	\leq -\gamma (\epsilon\rho) \|C(z)-C(\overline z)\|^2$, where we used the fact that $C^p$ is $\epsilon\rho$-cocoercive (cf. \autoref{prop:primal_sc}(iv)). 
	This implies that $\dot V(z)\! \leq 0,\! \forall z\in \mathbb{R}^n$ and $\dot V(z) \!= \!0$ only if $z \!\in \! \mathcal{E} \!:=\! \{z \in \mathcal{D}|C(z) \!=\! C(\overline z) \}$.  
	We find the largest invariant set $\mathcal{M}$ contained in $\mathcal{E}$ for $\dot z =\gamma (-z \!+\! U\circ C(z)).$ 
	On  $\mathcal{E}$, $\dot z \! 
	=\!  \gamma (-z +\overline z), $ hence since  $\gamma\!>0$,  $\|z(t) -  \overline z\|_\star \to 0$  as $ t\to \infty$, for any $z(0)\!\in \!\mathcal{E}$. Thus,  no other solution except $\overline z$ can stay forever in $\mathcal{E}$, and   $\mathcal{M}$ consists only of equilibria.  
	Since by assumption there are  a finite number of isolated equilibria, by LaSalle's invariance principle,  \cite{Khalil},
	it follows that for any $z(0)\! \in \!\mathcal{D}$, $z(t)$ converges to one of them,  $\overline z$. Finally,  since $C^p$ is $(\epsilon\rho)^{-1}$-Lipschitz (cf. \autoref{prop:primal_sc}(iii)), $ \|x(t) \!-\! \overline x\| \!\leq \!(\epsilon\rho)^{-1}\|z(t) \!-\!\overline z\|_\star $, hence $ \|x(t) - \overline x\| \to 0$ as $ t\to \infty$,  
	where,  by \autoref{prop:pertNash}, $ \overline x$  is a perturbed Nash equilibrium.  
			
	Alternatively, under \autoref{assump:pseudo_gradient}(iv), following from \eqref{eqn:strongly_convex_proof_1},\vspace{-0.2cm}
	\begin{align*}	
		\dot V(z) 	&\leq \gamma \mu\|C(z) - C(\overline z)\|^2 - \gamma ( C(z)\!-\! C(\overline z))^\top(z \!-\! \overline z)\\
					&\leq \gamma \mu\|C(z) - C(\overline z)\|^2  -\gamma (\epsilon\rho) \|C(z)-C(\overline z)\|^2\\
					&\leq -\gamma (\epsilon\rho - \mu)\|C(z) - C(\overline z)\|^2,
	\end{align*}
	where we again used the $\epsilon\rho$-cocoercivity of $C^p$. Assuming that $\epsilon\! >\! \mu\rho^{-1}$, then $\dot V(z) \!\leq \!0$, and convergence follows as  before. 
\end{proof}
\vspace{-0.3cm}

\begin{thm}
	\label{thm:convergence_Legendre}
	Let $\mathcal{G} \!= \!(\mathcal{N}, (\Omega^p)_{p \in \mathcal{N}}, (\mathcal{U}^p)_{p \! \in \!\mathcal{N}})$ be a concave game with players' dynamics given by DMD \eqref{eqn:DMD}. Assume there are a finite number of isolated fixed-points $\overline z$ of $U \circ C$, where $C \!=\! (C^p)_{p \in \mathcal{N}}$ is the mirror map induced by $\psi^p \!=\! \epsilon \vartheta^p$ satisfying \autoref{assump:primal}(ii). Then, under either \autoref{assump:pseudo_gradient}(i), (ii), or (iii),  with the additional assumption that $-\mathcal{U}^p$ is coercive in $x^p$ whenever $\Omega^p$ is non-compact, for any $\epsilon \!>\! 0$,  the auxiliary variables $z(t) \!=\! (z^p(t))_{p \in \mathcal{N}}$ converge to a rest point $\overline z$ while players' actions $x(t) \!=\! (x^p(t))_{p \in \mathcal{N}}$ converge to $\overline x$, a perturbed Nash equilibrium of $\mathcal{G}$. 
\end{thm}\vspace{-0.3cm}
\begin{proof}
	We use the same  Lyapunov function \eqref{eqn:global_convegence_bregman_sc}. 
	Since under  \autoref{assump:primal}(ii), ${\psi^p}^\star$ is Legendre (cf. \autoref{prop:primal_Legendre}(i)), ${\psi^p}^\star$ is strictly convex on $\interior(\dom {\psi^p}^\star)$, hence $V$ is positive definite at $z = \overline z$.  
Moreover, since ${\psi^p}^\star$ is essentially strictly convex, by \cite[Thm 3.7(iii)]{Bauschke97}, $D_{{\psi^p}^\star}(\cdot, \overline z^p)$ is coercive, so that  $V$ is radially unbounded. 
Then along any solution trajectory of  \eqref{eqn:DMD}, using   \autoref{prop:primal_Legendre}(ii), we can write as in \eqref{eqn:strongly_convex_proof_1} , 
$
	 \dot V(z)
	\!=\! \gamma (x \!-\! \overline x)^\top\!(u \!-\! \overline{u}) \!-\! \gamma (C(z) \!-\! C(\overline z))^\top\!(z \!-\! \overline z) $.  Since $u \!=\! U(x)$, $\overline u \! =\! U(\overline x)$, under either \autoref{assump:pseudo_gradient}(i), \ref{assump:pseudo_gradient}(ii) or \ref{assump:pseudo_gradient}(iii), the first term of $\dot V(z)$ is non-positive, so that, 
$	\dot V(z) \! \leq \! -\gamma(C(z) \!-\!C(\overline z))^\top\!(z\!-\!\overline z)$. 
	Since $C$ is strictly monotone  by \autoref{prop:primal_Legendre}(iii), therefore $\dot V(z)\! <\! 0, \forall z\! \in\! \mathbb{R}^n\backslash\{\overline z\}$, and by Lyapunov theorem \cite[Theorem 4.1, p.114]{Khalil}, $\overline z$ is asymptotically stable and therefore $z(t)$ converges to $\overline z, \!\forall z(0) \!\in \! \mathbb{R}^n$. By the continuity of $C(z)$ (\autoref{prop:primal_Legendre}(ii)), it follows that $x(t)$ converges  $ \overline x \!=\! C(\overline z), \forall x(0) \!=\! C(z(0))$, where  $\overline x$ is a perturbed Nash equilibrium. 
\end{proof}
\vspace{-0.3cm}
\begin{remark} In general, convergence is to the set of  perturbed Nash equilibria. By \autoref{prop:pertNash}, as $\epsilon \!\to \! 0$, $\overline x \! \to \! x^\star$, where $x^\star$ is a Nash equilibrium of $\mathcal{G}$. 
 Under \autoref{assump:pseudo_gradient}(ii) or \ref{assump:pseudo_gradient}(iii), the game admits a unique Nash equilibrium,  so $x(t)$ converges towards the unique $x^\star$. Note that in the case of Legendre regularizers, \autoref{thm:convergence_Legendre} gives convergence guarantees only for monotone games. 
 On the other hand, in the case of strongly convex regularizers, \autoref{thm:convergence_sc} gives guarantees for convergence in hypo-monotone games, 
 based on cocoercivity of the mirror map.
		We note that the above results can be extended to the weighted monotone case, \cite{Rosen65},
$			-(U(x) - U(x^\prime))^\top\Lambda(x-x^\prime) \geq 0, \forall x, x^\prime \in \Omega$,
		where $\Lambda = \diag(\lambda^1, \ldots, \lambda^N)$,  $\lambda^p > 0, p \in \mathcal{N}$, by appropriately redefining the regularizer. 
	\end{remark}

\vspace{-0.3cm}
\section{Examples of DMD}

We now provide several examples of DMD, whereby the mirror map is generated by regularizers in one of the two general classes. 
The first two are for examples of strongly convex regularizers (non-steep and steep), and the other three are for Legendre regularizers. For all derivations, we repeatedly use of the following result, based on a simple application of \cite[Theorem 4.14, p. 92]{Beck17}. 
\vspace{-0.2cm}
\begin{lem}

Let $\psi^p\! \!= \!\epsilon \vartheta^p\!$, $\epsilon \!>\!\! 0$, where  $\vartheta^p$ satisfies \autoref{assump:primal}, and let ${\psi^p}^\star$ be the convex conjugate of $\psi^p$. 
Then,    
	\begin{itemize} 
		\item[(i)] ${\psi^p}^\star(z^p) = \epsilon {\vartheta^p}^\star(\epsilon^{-1} z^p)$, 
		\item[(ii)] $C^p(z^p) = \nabla {\psi^p}^\star(z^p) = \nabla {\vartheta^p}^\star(\epsilon^{-1} z^p),  z^p \in \mathbb{R}^{n_p}$,
	\end{itemize} 
	where ${\vartheta^p}^\star$ is the convex conjugate of $\vartheta^p$. \vspace{-0.3cm}
\end{lem}


\vspace{-0.3cm}
\subsection*{Example 1. Euclidean Regularization over Compact Sets } 
 
 	Let $\Omega^p \subset \mathbb{R}^{n_p}$ be nonempty, compact and convex and consider, \vspace{-0.2cm}
	\begin{align} \label{eqn:primal_map_euclidean_norm}  & \vartheta^p(x^p) = \dfrac{1}{2}\|x^p\|^2_2 + \delta_{\Omega^p}(x^p). \end{align}
	By inspection, $\vartheta^p$  is supercoercive,  $1$-strongly convex  (\autoref{assump:primal}(i)) and non-steep, hence, $\psi^p \!=\! \epsilon \vartheta^p$  inherits the same properties  over $\Omega^p\!=\! \dom(\psi^p)$. 
	The convex conjugate  ${\psi^p}^\star$ 
	is given by, \vspace{-0.15cm}
	\begin{equation}
	\label{eqn:dual_phi_euclidean_projector}
	{\psi^p}^\star(z^p) = \dfrac{\epsilon}{2} (\|\epsilon^{-1} z^p\|_2^2 -  \|\epsilon^{-1} z^p  - \pi_{\Omega^p}(\epsilon^{-1} z^p)\|_2^2),
	\end{equation}
	where $\pi_{\Omega^p}$ is the Euclidean projection on $\Omega^p$. By \autoref{prop:primal_sc}(ii), ${\psi^p}^\star$ is continuously differentiable on $\mathbb{R}^{n_p}$ and can be shown to have a gradient $\nabla {\psi^p}^\star(z) = C^p(z^p) = \pi_{\Omega^p}(\epsilon^{-1} z^p)$. By \autoref{prop:primal_sc}(iii), (iv), (v), (vi), $C^p(z)$ is $\epsilon^{-1}$-Lipschitz, $\epsilon$-cocoercive, surjective from $\mathbb{R}^{n_p}$ onto $\Omega^p$ and has a left-inverse on $\Omega^p$ given by $\nabla \psi^p(x^p) = \epsilon x^p$. Then the DMD corresponding to \eqref{eqn:DMD},  \eqref{eqn:primal_map_euclidean_norm} is given by, \vspace{-0.2cm}
	\begin{equation}
	\label{eqn:PDMD}
	\begin{cases}
	\dot z^p &= \gamma (-z^p + u^p),\\
	x^p &= C^p(z^p) = \pi_{\Omega^p}(\epsilon^{-1} z^p),
	\end{cases}
	\end{equation}
	which we refer to as the \emph{projected DMD} (or PDMD).
	 By \autoref{thm:convergence_sc}, PDMD \eqref{eqn:PDMD} is guaranteed to converge to $\overline x$, a perturbed Nash equilibrium in any monotone game $\mathcal{G} = (\mathcal{N}, (\Omega^p)_{p \in \mathcal{N}}, (\mathcal{U}^p)_{p \in \mathcal{N}})$,  
	  for any $\epsilon > 0$,  and in any $\mu$ hypo-monotone game, for any $\epsilon > \mu$.

\begin{remark}  \eqref{eqn:PDMD} can be viewed as the continuous-time dual counterpart to the Tikhonov (primal) regularization algorithm in \cite{Shanbhag12}.
\end{remark}
\vspace{-0.3cm}
\subsection*{Example 2. Entropy Regularization over the Unit Simplex } 
Let $\Omega^p \! = \{x^p \! \in \!\mathbb{R}^{n_p}| x^p_i \geq 0, \sum\limits_{i = 1}^{n_p} x_i^p = 1\}:=\! \Updelta^p\!$   and 
\vspace{-0.23cm}
\begin{equation} \label{eqn:primal_map_simplex} \vartheta^p(x^p) =  \sum\limits_{i = 1}^{n_p} x_i^p \log(x_i^p).\end{equation}
with the convention $0\! \log 0 \!=\!0$. It can be shown that $\vartheta^p$ is supercoercive,  $1$-strongly convex over $\Updelta^p\!$ with respect to $\|\cdot\|_1$ (\autoref{assump:primal}(i))) and steep. Hence $\psi^p \!=\! \epsilon \vartheta^p$ inherits the same properties. 
Then 
${\psi^p}^\star(z^p) \!=\! \epsilon \log(\sum\limits_{i \!=\! 1}^{n_p}\! \exp(\epsilon^{-1} \!z^p))$, and  \vspace{-0.5cm} 
\begin{equation}
\label{eqn:EXPRL_map}
\nabla {\psi^p}^\star(z^p)\! =\! C^p(z^p) \!=\! \left [\dfrac{\exp(\epsilon^{-1} \!z^p_i)}{\sum\limits_{j \! = \!1}^{n_p} \exp(\epsilon^{-1} \!z_j^p)} \right ]_{i \in \{1,\dots,n_p\}}
\end{equation} 
By \autoref{prop:primal_sc}(iii), (iv), (v), (vi), $C^p(z^p)$ is $\epsilon^{-1}$-Lipschitz with respect to $\|\cdot\|_\infty$, $\epsilon$-cocoercive, surjective from $\mathbb{R}^{n_p}$ onto $\rinterior(\Updelta^p)$ and has a left-inverse on $\rinterior(\Updelta^p)$ given by $\nabla \psi^p(x^p) \!=\! \epsilon(\log(x^p) \!+\! \mathbf{1})$.  
Then the DMD corresponding to \eqref{eqn:primal_map_simplex} is given by  \eqref{eqn:DMD},  \eqref{eqn:EXPRL_map}, and  
	 by \autoref{thm:convergence_sc}, 
	is guaranteed to converge to a perturbed NE 
	  in any monotone game $\mathcal{G}\! =\! (\mathcal{N}, \!(\Updelta^p)_{p \in \mathcal{N}},\!(\mathcal{U}^p)_{p \in \mathcal{N}})$,   for any $\epsilon \!>\! 0$, and in any $\mu$ hypo-monotone game, for any $\epsilon \!>\! \mu$.
%
%
%
This dynamics corresponds to the exponentially-discounted reinforcement learning dynamics (EXP-D-RL) for finite games studied in \cite{Bo_LP_CDC2018}. 	There are several other well-known entropies over the simplex $\Updelta^p$ which are steep, e.g. 
	the log-barrier or the  Burg entropy, 
$		\vartheta^p(x) = - \sum\limits_{i = 1} \log(x_i^p), 
$ 
\cite{Mertikopoulos16}. \emph{Undiscounted} dynamics were shown to converge in games with a strict NE, \cite{Mertikopoulos16}, but not in zero-sum games with an interior NE. According to  \autoref{thm:convergence_sc}, the discounted DMD dynamics corresponding to these entropies are in fact guaranteed to converge in monotone games. 

\vspace{-0.3cm}
\subsection*{Example 3. Entropy Regularization over Non-Negative Orthant }   
Let $\Omega^p \!=\! \mathbb{R}^{n_p}_{\geq0}$ and consider the  \textit{Boltzmann-Shannon entropy} 
\vspace{-0.23cm}
	\begin{equation} 
	\label{eqn:primal_map_non_negative_orthant} 
	\hspace*{-0.1in} \vartheta^p(x^p) \!=\! \textstyle\sum\limits_{j =1}^{n_p} \left[ x^p_j \log(x^p_j) \!-\! x^p_j \right],
	\end{equation}
with the convention $0 \log 0 \!=\!0$. It can be shown that $\vartheta^p$ is supercoercive and Legendre\cite{Bauschke97} (\autoref{assump:primal}(ii)), hence $\psi^p \!= \! \epsilon \vartheta^p$ is too. 
	The dual map ${\psi^p}^\star\!: \!\mathbb{R}^{n_p} \to \mathbb{R} \!\cup \! \{\infty\}$ is given by, \vspace{-0.23cm}
	\begin{equation} 
	\label{eqn:dual_phi_sum_of_exponential} 
	{\psi^p}^\star(z^p) \!= \! \epsilon  \sum\limits_{i \!=\! 1}^{n_p} \exp(\epsilon^{-1}z^p_i), \, \text{ and } \, C^p(z^p) \!=\! \exp(\epsilon^{-1} z^p), 
	\end{equation}
which is  strictly monotone over $\interior(\mathbb{R}^{n_p}_{\geq 0})$. 
Then the DMD corresponding to \eqref{eqn:primal_map_non_negative_orthant} is  given by \eqref{eqn:DMD}, with $C^p$, \eqref{eqn:dual_phi_sum_of_exponential}, 
which we refer to as the \textit{Boltzmann-Shannon DMD} (or BDMD). 
By \autoref{thm:convergence_Legendre}, BDMD 
 is guaranteed to converge to $\overline x$, a perturbed Nash equilibrium in any monotone game $\mathcal{G} \!=\! (\mathcal{N}, \!(\mathbb{R}^{n_p}_{\!\geq 0})_{p \! \in \!\mathcal{N}}, (\mathcal{U}^p)_{p \!\in\! \mathcal{N}})$,  
	  for any $\epsilon \!>\! 0$. 
	  %
The BDMD can be generalized to $\Omega^p \!=\! [-c^p, \infty]^{n_p}, c^p \!\geq \!0$ 
	and the mirror map is given by $C^p(z^p)\! =\! \exp(\epsilon^{-1}\!z^p) \!-\! c^p\mathbf{1}$. 

\vspace{-0.3cm}
\subsection*{Example 4.  Entropy Regularization over Unit Square }   

Let $\Omega^p = [0,1]^{n_p}$ and consider the \textit{Fermi-dirac entropy} \vspace{-0.23cm} 
\begin{equation} 
\label{eqn:primal_map_unit_square}
\vartheta^p(x^p) = \textstyle\sum\limits_{j =1}^{n_p} \left[ (x^p_j ) \log(x^p_j ) \!+\! (1-x^p_j)\log(1-x_j^p) \right],
\end{equation} 
which can be shown to be supercoercive and Legendre \cite{Bauschke97}  (\autoref{assump:primal}(ii)), hence $\psi^p \!=\! \epsilon \vartheta^p$ is supercoercive and Legendre as well. The dual map ${\psi^p}^\star: \mathbb{R}^{n_p} \to \mathbb{R}\cup\{\infty\}$ is given by, \vspace{-0.23cm}
\begin{equation} {\psi^p}^\star(z^p) =   \epsilon \sum\limits_{i = 1}^{n_p}  \log(1+\exp(\epsilon^{-1} z_i^p)), \end{equation} 
sometimes referred to as the \textit{softplus} function. The mirror map is \vspace{-0.23cm}   
\begin{equation}\label{eqn:FDMD_map}
C^p(z^p) = \left[\dfrac{\exp(\epsilon^{-1} z^p_i)}{1+\exp(\epsilon^{-1} z^p_i)} \right]_{i \in \{1,\dots,n_p\}}
\end{equation}
It can be shown that $C^p$ is strictly monotone over $\interior([0, 1]^{n_p})$  with inverse $\nabla \psi^p(x^p) \!=\! (\epsilon \log(x_j^p/(1-x_j^p)))_{j = 1}^{n_p}.$ The associated DMD  
 is given by \eqref{eqn:DMD}, \eqref{eqn:FDMD_map},  
which we refer to as the \textit{Fermi-Dirac regularized} (FDMD). 
By \autoref{thm:convergence_Legendre}, FDMD \eqref{eqn:DMD}, \eqref{eqn:FDMD_map},   is guaranteed to converge to a perturbed NE in any monotone game $\mathcal{G} \!=\! (\mathcal{N}, ([0,1]^{n_p})_{p \in \mathcal{N}}, (\mathcal{U}^p)_{p \in \mathcal{N}})$,  
	  for any $\epsilon \!>\! 0$. 
The FDMD can be generalized to $\Omega^p \!=\! [a^p, \!b^p]^{n_p}$, by appropriately modifying it. 
%


\begin{table*}[ht!]
	\centering
	\begin{tabular}{Lcllc}
		\toprule   
		\textbf{Name and Acronym} &   \textbf{Dynamics} & \textbf{Mirror Map} & \textbf{Player Action Set}  \\   
		\toprule
		\midrule
		Projected Discounted MD (PDMD)
		& $\dot z^p = \gamma(-z^p + u^p)$ 
		& $ x^p =  \pi_{\Omega^p}(\epsilon^{-1} z^p)$ & $\Omega^p$\\  
		Exponentially-Discounted RL (EXPD-RL)
		& 
		& 
		$ x^p =   \left [ \dfrac{\exp(\epsilon^{-1} z^p_i)}{\sum\limits_{j = 1}^{n_p} \exp(\epsilon^{-1} z_j^p)} \right ]_{i \in \{1,\dots,n_p\}}$ & $\Delta^p$\\ 
			Boltzmann-Shannon Regularized DMD (BDMD)
		& 
		&  $x^p = \exp(\epsilon^{-1}z^p) - c^p\mathbf{1}$ & $[-c^p, \infty]^{n_p}$ \\
		Fermi-Dirac Regularized DMD (FDMD) 
		& 
		&
		  $x^p =  \left[\dfrac{a^p\!+\!b^p\exp(\epsilon^{-1} z_i^p) }{\exp( \epsilon^{-1} z_i^p) + 1}\right]_{i \in \{1,\dots,n_p\}}$ & $[a^p, b^p]^{n_p}$ \\ 
		Hellinger Regularized DMD (HDMD) 
		& 
		& $x^p  = a^p\epsilon^{-1} z^p\left[\sqrt{1+\|\epsilon^{-1}z^p\|_2^2}\right]^{-1} - c^p\mathbf{1} $ & $\overline{\mathcal{B}}^{n_p}_{a^p}(c^p)$\\
		\bottomrule
	\end{tabular}
	\caption{Discounted Mirror Descent Dynamics} 
	\label{table:list_md_dynamics}
\end{table*}
\vspace{-0.3cm}
\subsection*{Example 5. Regularization over Euclidean Spheres}  

Assume that $\Omega^p \! = \! \{x^p \!\in \!\mathbb{R}^{n_p}| \|x^p \!-\! c^p\|_2 \leq  a^p\}\! :=\!\overline{\mathcal{B}}^{n_p}_{a^p}(c^p)$.  
Consider the \textit{Hellinger distance}, 
$\vartheta^p(x^p) = -\sqrt{{a^p}^2 - \|x^p - c^p\|_2^2},
$ 
which can be shown to be supercoercive and Legendre\cite{Bauschke97}. Hence $\psi^p \!=\! \epsilon \vartheta^p$ is supercoercive and Legendre as well  (\autoref{assump:primal}(ii)). The dual map is 
${\psi^p}^\star(z^p) = \epsilon (a^p\sqrt{1 + \|\epsilon^{-1} z^p\|_2^2} - \epsilon^{-1}{c^p}^\top z)$ and  
the mirror map is given by, \vspace{-0.23cm}
\begin{equation}\label{eqn:HDMD_map}
C^p(z^p) = \dfrac{a^p \epsilon^{-1} z^p}{\sqrt{1+\|\epsilon^{-1}z^p\|_2^2}} - c^p\mathbf{1},
\end{equation}
 strictly monotone over $\mathbb{R}^{n_p}$. 
The associated DMD 
 given by \eqref{eqn:DMD}, \eqref{eqn:HDMD_map}, 
which we refer to as the \textit{Hellinger regularized} DMD (HDMD),   
is guaranteed to converge to a perturbed NE in any monotone game $\mathcal{G} \!=\! (\mathcal{N},  (\overline{\mathcal{B}}^{n_p}_{a^p}(c^p))_{p \in \mathcal{N}}, (\mathcal{U}^p)_{p \in \mathcal{N}})$,  
	  for any $\epsilon \!>\! 0$, cf.  \autoref{thm:convergence_Legendre}.
	  
	We summarize all these discounted dynamics in \autoref{table:list_md_dynamics}. Note that the undiscounted versions of these dynamics are given by \eqref{eqn:mirror_descent} with the corresponding mirror maps. 
	

\section{Simulation Results}
In this section, we provide simulations results.  We note that an example of resource sharing via Kelly's mechanism \cite{Kelly98} in a 
strictly monotone game for $N\!=\!4$ players is provided in \cite{Bo_LP_CDC2019}. Here consider representative examples of monotone and hypo-monotone games. For comparison purposes, all dynamics are simulated over the same duration, and  
unless otherwise specified,  
with the same $\epsilon$ and  $\gamma \!\!=\!\! 1$, for initial value $z(0) \!\!=\!\! \mathbf{0}$. For PDMD \eqref{eqn:PDMD} and FDMD 
 we assume that each player's strategy is projected onto $[-100, \!100]$ (or Cartesian product of it) whenever the action set is unconstrained, and onto  $[0,\! 100]$ whenever the action set is a subset of the 
non-negative orthant. For HDMD, 
we used 
the ball of radius $100$ centered at  origin. The color code for each dynamics is as follows: PDMD \eqref{eqn:PDMD} (blue), BDMD 
(red), FDMD  
(orange), HDMD 
(magenta). 

\begin{figure}[ht]
\centering
\begin{minipage}[ht]{0.45\columnwidth}
\hspace{-0.2cm}	\centerline{\includegraphics[width=1.15\linewidth]{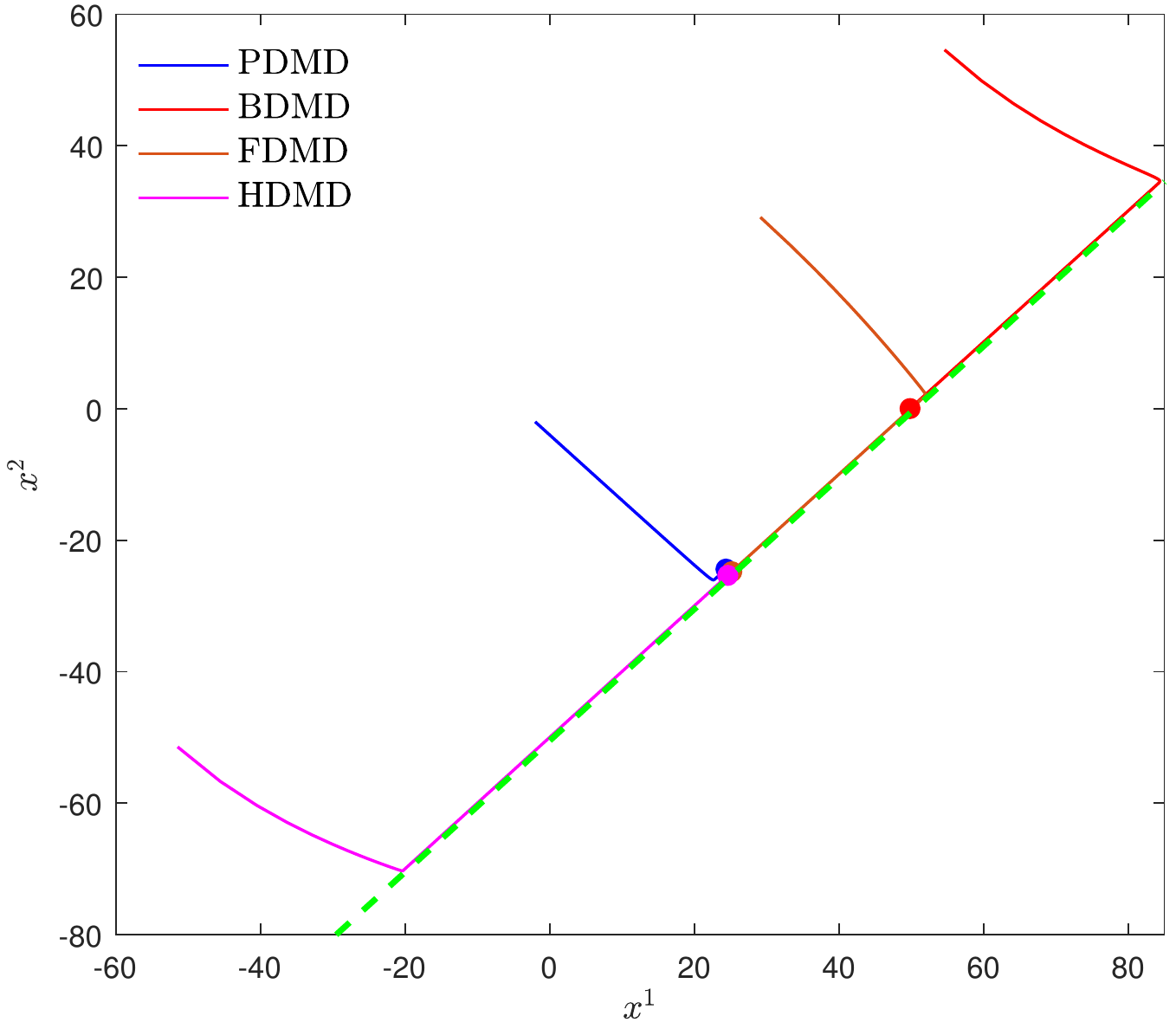}}
	\caption{Monotone game, $\epsilon\!\!=\!\!0.5$}  
	\label{fig:quadratic_game_eps_0_5}
	\end{minipage}
\hspace{0.25cm}
\begin{minipage}[ht]{0.45\columnwidth}
\vspace{0.27cm}	
\centerline{\includegraphics[width=1.15\linewidth]{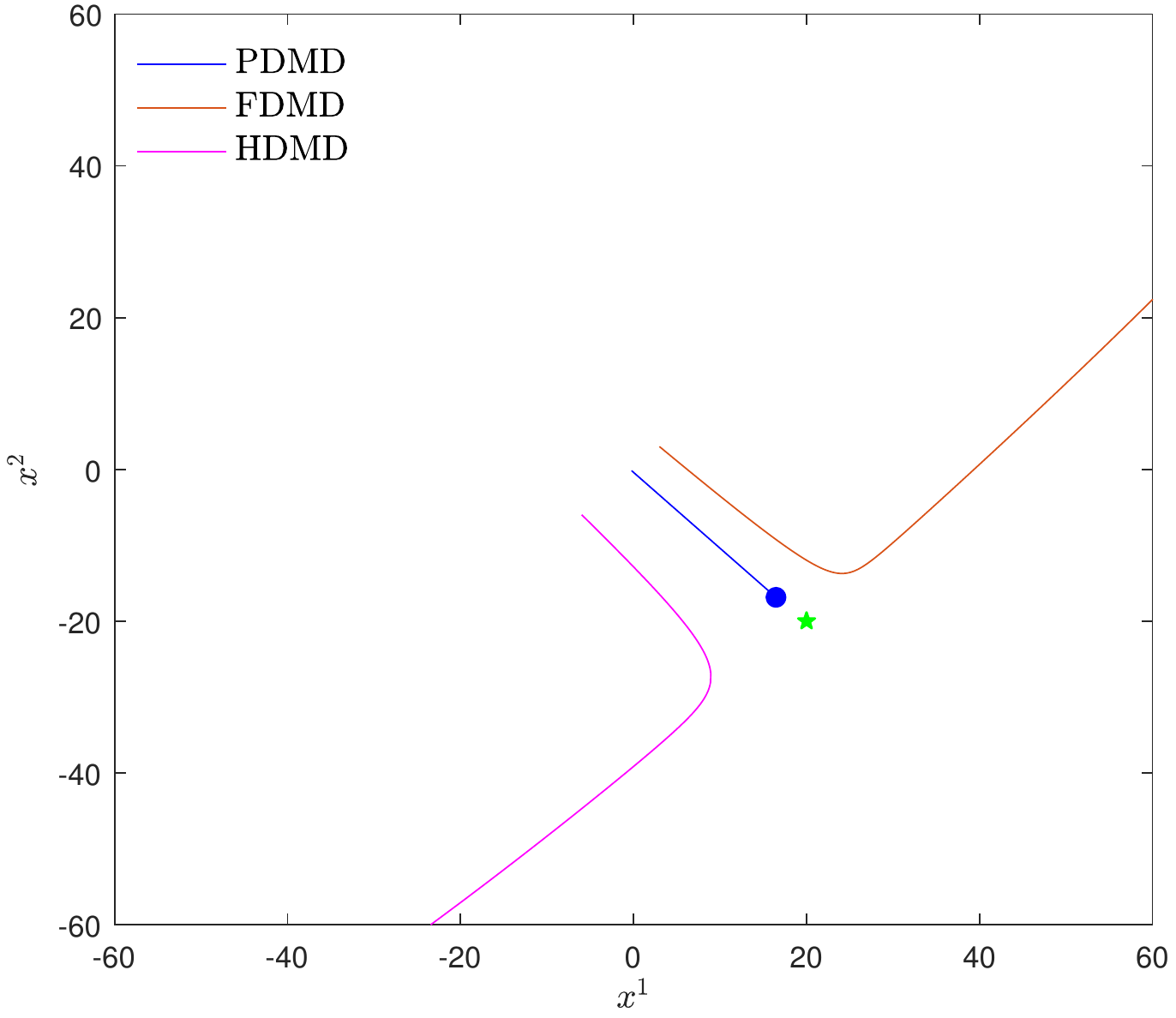}}
	\caption{Hypo-monotone game, $\epsilon\!\!=\!\!5.1$}     
	\label{fig:quadratic_game_eps_5_hypo}
\end{minipage}
\end{figure}

\begin{example}(\textbf{Monotone and Hypo-monotone Quadratic Games})\\
	In this example we compare the discounted DMDs in a monotone quadratic and a hypo-monotone quadratic game. For the monotone game, we also compare them with the discrete-one introduced in \cite{Shanbhag10, Shanbhag12}.  
	Quadratic games constitute an important class of games that serve as second-order approximation to other nonlinear payoff functions and models of competition between markets \cite[p. 190]{Basar}. Consider an $N$ player game where each player $p$ has a quadratic payoff function, 
$	\mathcal{U}^p(x^p; x^{-p}) \!\!=\!\! \dfrac{1}{2} \sum\limits_{i = 1}^N \sum\limits_{j = 1}^N {x^i}^\top A^p_{ij} x^{j} \!\!+\!\! \sum\limits_{i = 1}^N b^p_i x^i \!\!+\!\! c^p
$, 
	where   $x^p \in \mathbb{R}^{n_p}$ and  $A_{ij}^p \!\!=\!\! {A_{ji}^p}^\top \!\! \in \!\! \mathbb{R}^{n_p\times n_p}$, with each $A^p_{ii}$ being symmetric and $b^p_i  \!\!\in \!\! \mathbb{R}^{n_p}$, $c^p \!\!\in \!\! \mathbb{R}$. Let $A^p \in \mathbb{R}^{n \times n}$ be the block matrix, $A^p \!\!=\!\! (A^p_{ij})$, and $b^p \!\!=\!\! [ b^p_1 \ldots b^p_N ]^\top \!\! \in \!\! \mathbb{R}^n$, we can write 
$	\mathcal{U}^p(x^p; x^{-p}) \!=\!\dfrac{1}{2} x^\top A^p x \!+\!{b^p}^\top x \!+\! c^p. 
$
The pseudo-gradient  of this game is \vspace{-0.3cm}
	\begin{equation} \label{eqn:quadratic_game_pseudo_gradient} U(x) = Rx + b, \end{equation} 	where $R = \begin{bmatrix} A^1_{11} & \ldots & A^1_{1N} \\ \vdots & \ddots & \vdots \\ A^N_{N1} & \ldots & A^N_{NN} \end{bmatrix}$. 
	Then the game is monotone (cf. \autoref{assump:pseudo_gradient}(i) if for all  $x, x^\prime \! \in\!  \mathbb{R}^{n}$, 
$	- (U(x) \!-\! U(x^\prime))^\top(x-x^\prime)  \!=\!  -(x\! -\! x^\prime)^\top R (x-x^\prime)\! \geq \!0 
$, i.e., if 
	$R \!+\!R^\top$ is negative semidefinite.
	
	Consider $N \!=\! 2$, $x^p \!\in \! \mathbb{R}$, 
	$A^1 \!\! = \!\! \begin{bmatrix} -10 & \!10 \\ 5 & \! -5 \end{bmatrix}$, $A^2 \!\!=\!\! \begin{bmatrix} -5 & \! 5 \\ 10  & \! -10 \end{bmatrix}$,  $b^1 \!\!= \!\! 500, b^2 \!\!=\!\! -500$, $c^1 \!\!=\!\! 0$, $c^2 \!\! = \!\! 0$. 
	The pseudo-gradient $U(x)$ is \eqref{eqn:quadratic_game_pseudo_gradient} where $R \!\!=\!\! \begin{bmatrix}  -10 & \! 10 \\ 10 & \!-10 \end{bmatrix}$, $b  \!\!=\!\! \begin{bmatrix} 500 & -500 \end{bmatrix}^\top $. 
	$R$ has eigenvalues $\{-20, 0\}$, 
	hence the game is monotone. The set of Nash equilibria is $x^\star \!\in \! \{(x^1, \!x^2) \in \mathbb{R}^2| x^1 \!\!=\!\! 50 \!+\! x^2\}$, set indicated with a green line in \autoref{fig:quadratic_game_eps_0_5}. 
	In 	\autoref{fig:quadratic_game_eps_0_5}, we provide simulations of  PDMD, 
	BDMD, 
	FDMD 
	and HDMD,   
	all for  $\epsilon = 0.5$, in $(x^1, x^2)$ plane. 
	In order to distinguish between trajectories, each dynamics is simulated with a different initial $z(0)$.  We observe that each of the dynamics PDMD, FDMD and HDMD converges close to  
	$(25,\!-25)$ (a NE), while  BDMD converges close to 
	$(50, \!0)$ (also a NE). 
	



\end{example}

\vspace{-0.5cm}	
\vspace{-0.27cm}
Consider now 
$A^1 \! \!= \!\! \begin{bmatrix} -10 & \!15 \\ 15 & \! -5 \end{bmatrix}$, $A^2 \!\!= \!\!\begin{bmatrix} -5 & \!15 \\ 15  & \!-10 \end{bmatrix}$, so that 
the pseudo-gradient, \eqref{eqn:quadratic_game_pseudo_gradient}, has 
$R \!\!=\!\! \begin{bmatrix}  -10 & \!15 \\ 15 & \!-10 \end{bmatrix}$, 
with eigenvalues $\{-25, 5\}$. By  \autoref{assump:pseudo_gradient}(iv), the game is hypo-monotone ($\mu \! = \! 5$), and 
$x^\star \! = \! (20,\! -20)$. 
We note that only the PDMD  is guaranteed to converge to a perturbed NE (cf. \autoref{thm:convergence_sc}), for $\epsilon \!\!>\!\! 5$. 
In \autoref{fig:quadratic_game_eps_5_hypo} we provide simulations of the DMDs 
for $\epsilon \!\! =\!\! 5.1$, which show that PDMD  converges to $(\overline{x}^1, \overline{x}^2) \!\! = \!\!(16.4,\! -16.7)$ as per  \autoref{thm:convergence_sc} (relatively close to $x^\star$ (green star)), while the other dynamics fail to converge. 

\vspace{-0.1cm}
\begin{remark}
	We compare the discretization 
	of PDMD \!\eqref{eqn:PDMD} to 
	the (coordinated) iterative Tikhonov regularization (ITR) scheme \cite{Shanbhag12},  
	\begin{equation}
		x^p_{k+1} = \pi_{\Omega^p}(x^p_{k} - t_k (-U^p(x_k) + \epsilon_k x^p_k)),  
	\end{equation}
	where, $t_k$ and $\epsilon_k$ are sequences of diminishing step-size. By \cite[Theorem 2]{Shanbhag10}, ITR converge to the least-norm solution of $\text{VI}(\Omega, -U)$  (in the sense of \cite[p. 1128]{Facchinei_I})  for monotone games ($-U$ is monotone) with Lipschitz pseudo-gradient map when the $t_k, \epsilon_k$ are appropriately chosen \cite[Lemma 4]{Shanbhag10}. 	
	We compare ITR to the discrete-time PDMD obtained by an Euler discretization of \eqref{eqn:PDMD}, \vspace{-0.2cm}  
	\begin{equation} 
	\label{eqn:PDMD_discrete}
	 			\begin{cases} z_{k+1}^p &= z_k^p +  t_k (-z_k^p + U^p(x_k) ), \\     
	 			x_{k+1}^p &=   \pi_{\Omega^p}(\epsilon^{-1} z^p_{k+1}), \end{cases}
	\end{equation}
	where we use $t_k = 0.001$ and  $\epsilon = 0.1$. 
 	We run this for the monotone game considered before and we also run ITR with  $t_k = k^{-0.48},  \epsilon_k = k^{-0.51}$ (as in \cite{Shanbhag12}). 
The evolution of $x^p_k$ under PDMD and ITR is shown in \autoref{fig:quadratic_game_Shanbhag_monotone_comparison_constant},  
with ITR shown in orange and PDMD shown in blue (solid line one player, dashed line the other one). Both converge close to the NE at $x^\star = (25, -25)$, 
but we  find for all the step-sizes, PDMD has a faster rate of convergence  as compared to ITR. 

%
%
	\vspace{-0.4cm}
	\begin{figure}[htp!]
	\centering
	\includegraphics[width= 0.9\linewidth]{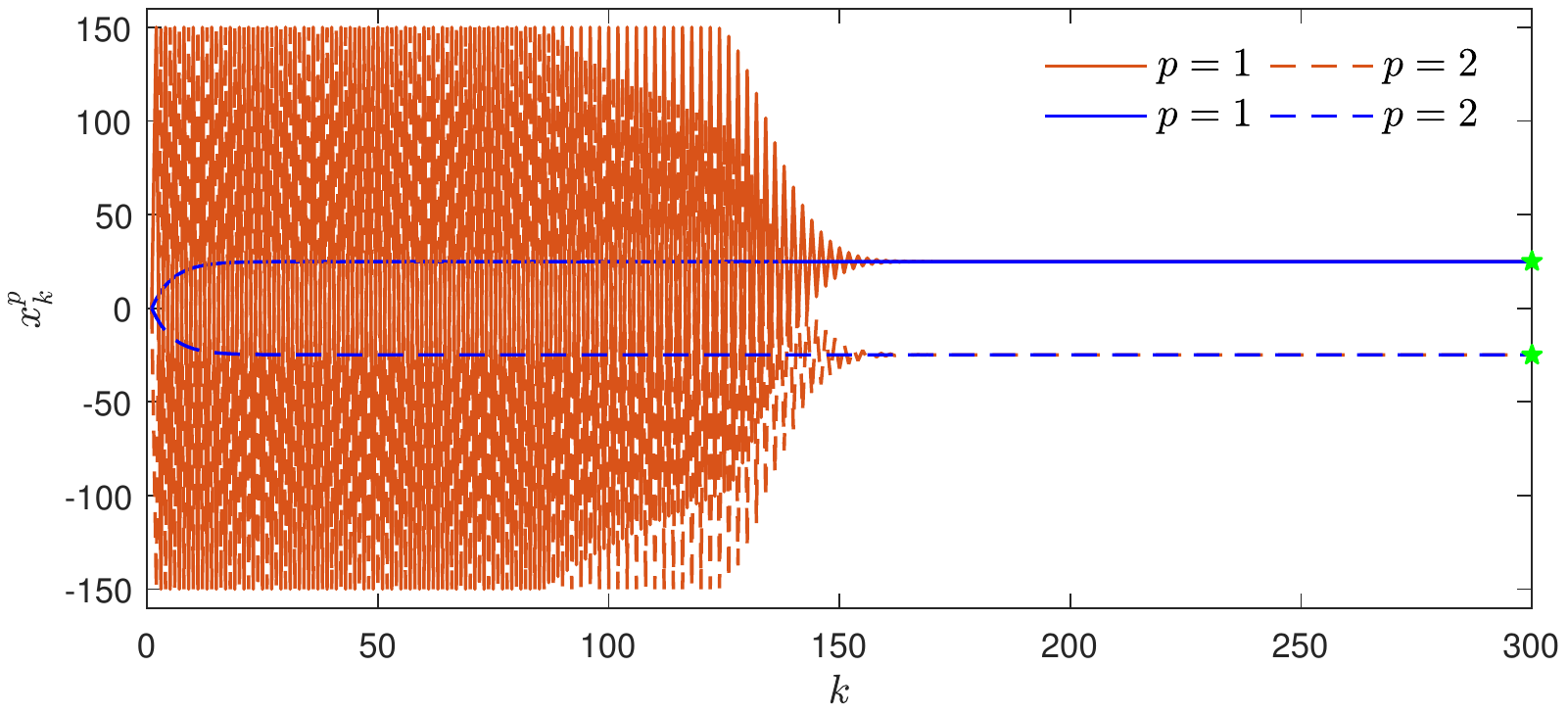}
	\caption{Comparison between discrete-time PDMD and ITR}   
	\label{fig:quadratic_game_Shanbhag_monotone_comparison_constant}
	\vspace{-0.2cm}
\end{figure}

	
\end{remark}

	\vspace{-0.25cm}
\begin{example}(\textbf{Learning the Mean of a Distribution})\\
In this example we compare discounted DMDs with their undiscounted versions, applied to learning the mean of a distribution, formulated as a monotone game. 
	Let $\mathsf{Z}\! \!\sim \!\!\mathsf{P}(\mathsf{z})$ and $\mathsf{X} \! \!\sim \!\!\mathsf{Q}(\mathsf{x}) $ be two random variables. We wish to construct a function $G_\theta\!: \!\mathbb{R}^n \! \to \!\mathbb{R}^n$, parameterized by an unknown parameter $\theta \! \in \! \mathbb{R}^n$, such that $\mathbb{E}(G_\theta(\mathsf{Z})) \!=\! \mathbb{E}(\mathsf{X})$. The authors of \cite{Daskalakis18} showed that $G_\theta$ can be constructed by solving the saddle point problem
	\footnote{
	This is an example of the so-called Generative Adversarial Network (GAN), specifically, the Wasserstein GAN (without Lipschitz constraint).
	}, 
	\vspace{-0.2cm}
	\begin{equation}
	\underset{\theta \in \mathbb{R}^n}{\text{min}} \thinspace \underset{w \in \mathbb{R}^n}{\text{max}} \thinspace \mathbb{E}_{\mathsf{X} \sim \mathsf{Q}(\mathsf{x})}(D_w(\mathsf{X})) - \mathbb{E}_{\mathsf{Z} \sim \mathsf{P}(\mathsf{z})}(D_w(G_\theta(\mathsf{Z})))  
	\label{eqn:gan_objective}
	\end{equation} 
	where $D_w \!:\! \! \mathbb{R}^n \!\! \to \!\!\mathbb{R}$ is a function parametrized by unknown parameter $w\! \!\in\!\! \mathbb{R}^n$. 
As an example, let $\mathsf{Z}  \!\!\sim\! \!\mathcal{N}(0,1)$, $\mathsf{X} \!\!\sim \!\!\mathcal{N}(v, 1)$ (Gaussian distributions), $v \!\!= \!\!\mathbb{E}(\mathsf{X}) \!\!\in \!\! \mathbb{R}$, $D_w(\mathsf{X}) \!\!= \!w^\top \mathsf{X}$ and $G_\theta(\mathsf{Z}) \!\!=\!\! \mathsf{Z}+\theta$, for $\theta, w  \!\in \!\! \mathbb{R}$. Then the objective in \eqref{eqn:gan_objective} is given by, 
			$ \mathbb{E}_{\mathsf{X} }(D_w(\mathsf{X})) \! -\! \mathbb{E}_{\mathsf{Z}}(D_w(G_\theta(\mathsf{Z})))
			\!=\! \thinspace w^\top \mathbb{E}_{\mathsf{X} }(\mathsf{X}) \!-\! w^\top \mathbb{E}_{\mathsf{Z} }(\mathsf{Z} \!+\! \theta) 
			\!=\! w^\top(v \!-\! \theta)
			$. 
    With $x^1 \!=\! \theta, x^2  \!= \!w$, \eqref{eqn:gan_objective} is equivalent to a two-player zero-sum game with \vspace{-0.2cm}
	\begin{equation}
	\mathcal{U}^1(x^1; x^2) = -{x^2}^\top(v - x^1) \qquad \mathcal{U}^2(x^2; x^1) = {x^2}^\top(v - x^1)
	\end{equation}
	where the player sets are $\Omega^1 = \mathbb{R}, \Omega^2 = \mathbb{R}$. The pseudo-gradient is
$	U(x) \!=\!R \begin{bmatrix} x^1 \\ x^2 \end{bmatrix}+\begin{bmatrix} 0  \\ v \end{bmatrix}$, $R\!=\! \begin{bmatrix} 0 & 1 \\ -1 & 0 \end{bmatrix}$, 
	hence 
	the game is monotone, and has NE $x^\star  \!=\! (v, \!0)$. 
Let $v\!=\!50$. 
	In \autoref{fig:Mean_eps_0.1}, we show results for the discounted dynamics  (solid),  
	as well as for their undiscounted counterparts (dashed), 
	for $\epsilon \!=\! 0.1$.
		We slightly increased the solver step-size for  FDMD  and HDMD in order to distinguish trajectories. 
As seen, 
all discounted DMD dynamics converge to the NE $x^\star$ (shown by a green star), whereas the undiscounted dynamics cycle.  

	

	\vspace{-0.15cm}	
	\begin{figure}[htp!]
		\centering
		\includegraphics[width= 0.9\linewidth]{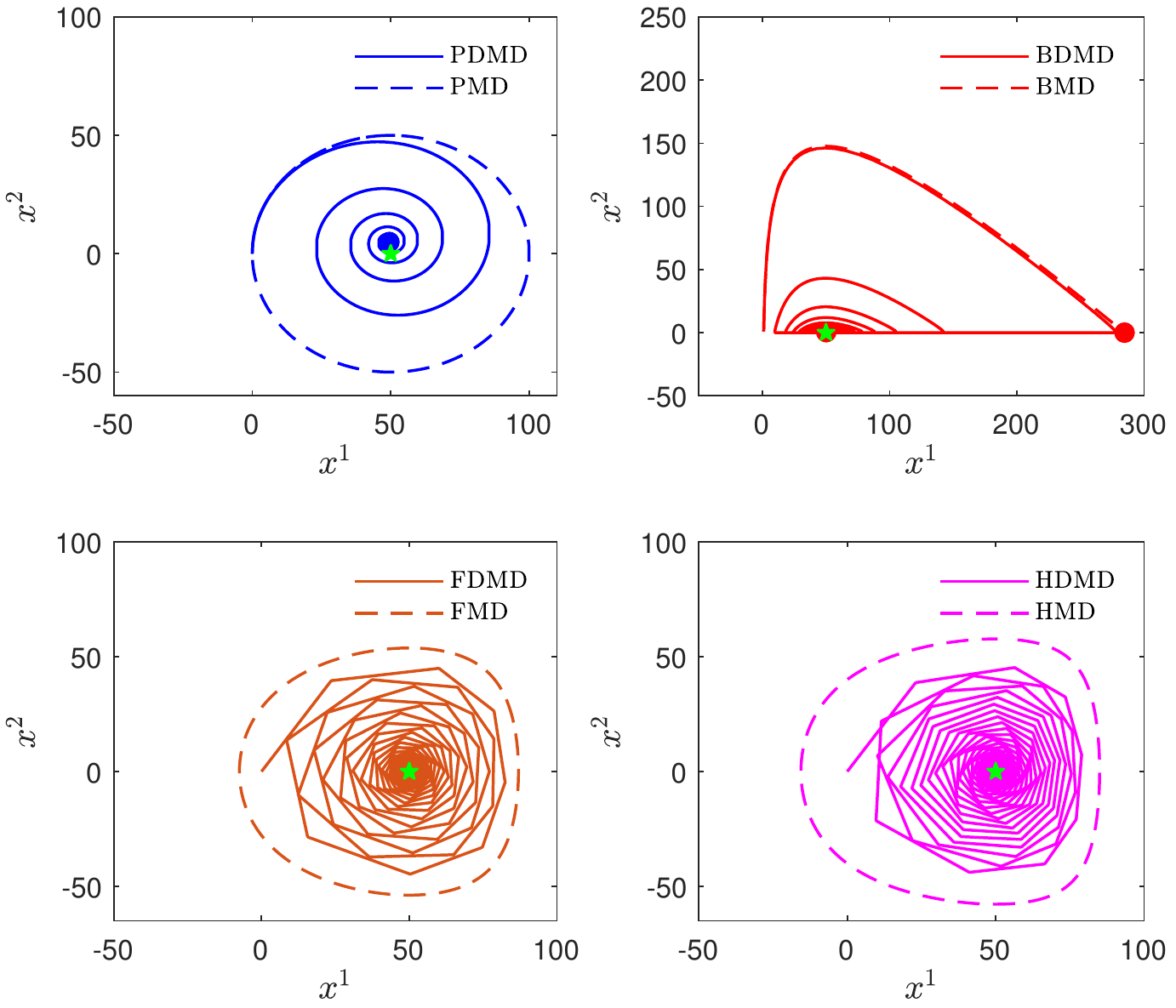}
		\caption{Learning the Mean, $\epsilon = 0.1$}
		\label{fig:Mean_eps_0.1}
		\vspace{-0.4cm}
	\end{figure}
\end{example}

\begin{example}(\textbf{Polynomial Regression})\\
In this example we compare the various DMD dynamics for a polynomial regression problem formulated as a zero-sum monotone game. 
	Consider a data set $\mathsf{D} \!=\! \{(a_n, b_n)\}_{n = 1}^\mathsf{N}$, $a_n \in \mathbb{R}, b_n \in \mathbb{R}$ and 
	\vspace{-0.25cm}
	\begin{equation}
	\label{eqn:polynomial}
	p(a) = \mathsf{w}_0 + \mathsf{w}_1 a + \mathsf{w}_2 a^2 + \ldots \mathsf{w}_\mathsf{M} a^\mathsf{M}, 
	\end{equation} where  the coefficients $\mathsf{w} \!=\! (\mathsf{w}_0, \mathsf{w}_1, \ldots, \mathsf{w}_\mathsf{M}) \in \mathbb{R}^{\mathsf{M}+1}$ 
are to be found for the best $M$-order  fit through the data $\mathsf{D} $,	$a \in \{a_1, \ldots, a_\mathsf{N}\}$. 
 These coefficients can be found by solving $\min\limits_{\mathsf{w} \in \mathbb{R}^{M+1}} \dfrac{1}{2}\|A\mathsf{w} - b\|_2^2 $, where  $A \!\in \!\mathbb{R}^{\mathsf{N} \times {\mathsf{M}+1}}$, 
%
$		A \!=\! \begin{bmatrix} 1 & a_1 & a_1^2 & \ldots & a_1^\mathsf{M} \\ \vdots &  \vdots  &  \vdots  &  \vdots  &  \\ 1 & a_\mathsf{N} & a_\mathsf{N}^2 & \ldots &  a_\mathsf{N}^\mathsf{M} \end{bmatrix}, 
		b \!=\! \begin{bmatrix} b_1 \\ \vdots \\ b_\mathsf{N} \end{bmatrix}.
$ 
		Assume $\range(A) \!=\! \mathbb{R}^\mathsf{N}$, 
		then the objective function can be rewritten as 
$	\dfrac{1}{2}\| A\mathsf{w} - b\|_2^2 \!=\! \max\limits_{\mathsf{z} \in \mathbb{R}^\mathsf{N}} \mathsf{z}^\top( A\mathsf{w} \!-\! b) \!-\! \dfrac{1}{2}\|\mathsf{z}\|_2^2
\!:=\! \! \max\limits_{\mathsf{z} \in \mathbb{R}^\mathsf{N}} f(\mathsf{w},\mathsf{z})$, \cite{Mokhtari19, Du19}, 
and 
$	\min\limits_{\mathsf{w} \in \mathbb{R}^{\mathsf{M}+1}} \dfrac{1}{2}\|A\mathsf{w} - b\|_2^2 = \min\limits_{\mathsf{w} \in \mathbb{R}^{M+1}} \max\limits_{\mathsf{z} \in \mathbb{R}^{N}}  f(\mathsf{w},\mathsf{z}) 
$. 
	With  $x^1 \!=\! \mathsf{w}\! \in\! \mathbb{R}^{\mathsf{M}+1}, x^2 \!=\! \mathsf{z} \!\in \! \mathbb{R}^\mathsf{N}$, we obtain a two-player zero-sum game with 
	payoff functions 
$	\mathcal{U}^1(x^1; x^2) \!=\! -f(x^1, x^2)$, $\mathcal{U}^2(x^2; x^1) \!=\! f(x^1, x^2)$.  
	The pseudo-gradient is 
$	U(x) \!=\! R \begin{bmatrix} x^1 \\ x^2 \end{bmatrix} \!-\! \begin{bmatrix} 0 \\ b \end{bmatrix} $, for $R\!=\!  \!\begin{bmatrix} \varnothing & \!-A^\top \\ A &  \!-\! I \end{bmatrix} $. 
Since $R+R^\top =\! \begin{bmatrix} \varnothing & \!\varnothing \\ \varnothing & \!-2I \end{bmatrix}$	is negative semidefinite,  the game is monotone. 
$U(x^\star) \!= \!0$ yields 
$		x^{\star 2} \! \in \! \ker(A^\top)$, 
$		Ax^{\star 1} \!-\! x^{\star 2}  \!-\! b\! =\! 0
$,
	hence $x^{\star 1} \!=\! (A^\top A)^{-1} A^\top b$. 
	Consider a data set with $\mathsf{N} \!=\! 20$  points, $\mathsf{D} \!=\! \{ (1,20), \!(2, 12), \!(3, 15), \ldots, \!(18, -10), \!(19,20), \!(20, 2)\}$ and $\mathsf{M} \!= \!3$ (fit to a third-order polynomial). The optimal coefficients are $x^{\star1}\! =\!  (36.0640, \!-13.0372, \!1.2084, \!-0.0342)$. 
	 \autoref{fig:Polynomial_Regression_Epsilon_0_1_Third_Order} shows $x^1$ trajectories under PDMD, FDMD  and HDMD, all with $\epsilon \! =\! 0.1$, as well as PSGD, 
	 with $ x^{\star1}$ as green stars. 
	We see that FDMD  and HDMD converge to $x^\star$,  while  PDMD 
	is very slow. The third-order polynomial associated with the final coefficients found by each dynamics, along with the best fit 
	are shown in \autoref{fig:Polynomial_Regression_Curve_Fit_Third_Order_Epsilon_0_1} (data points as red circles), indicating superior performance of  FDMD and HDMD. 

\vspace{-0.23cm}
\begin{figure}[ht]
\centering
\begin{minipage}[ht]{0.45\columnwidth}
\hspace{-0.2cm}	\centerline{\includegraphics[width=1.15\linewidth]{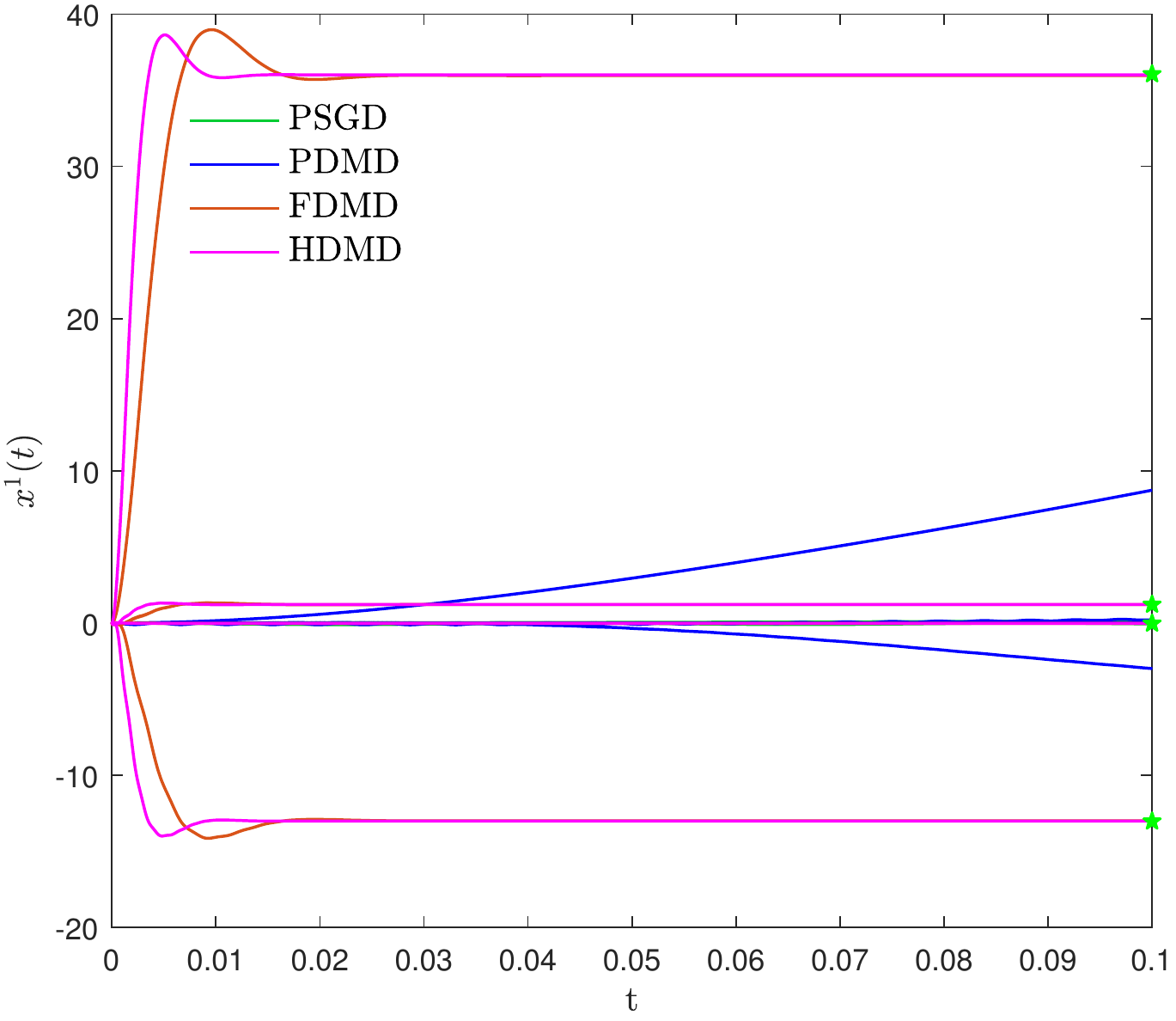}}
\vspace{-0.2cm}	\caption{Solution trajectories}  
	\label{fig:Polynomial_Regression_Epsilon_0_1_Third_Order}
	\end{minipage}
\hspace{0.2cm}
\begin{minipage}[ht]{0.45\columnwidth}
\vspace{-0.1cm}	\centerline{\includegraphics[width=1.1\linewidth]{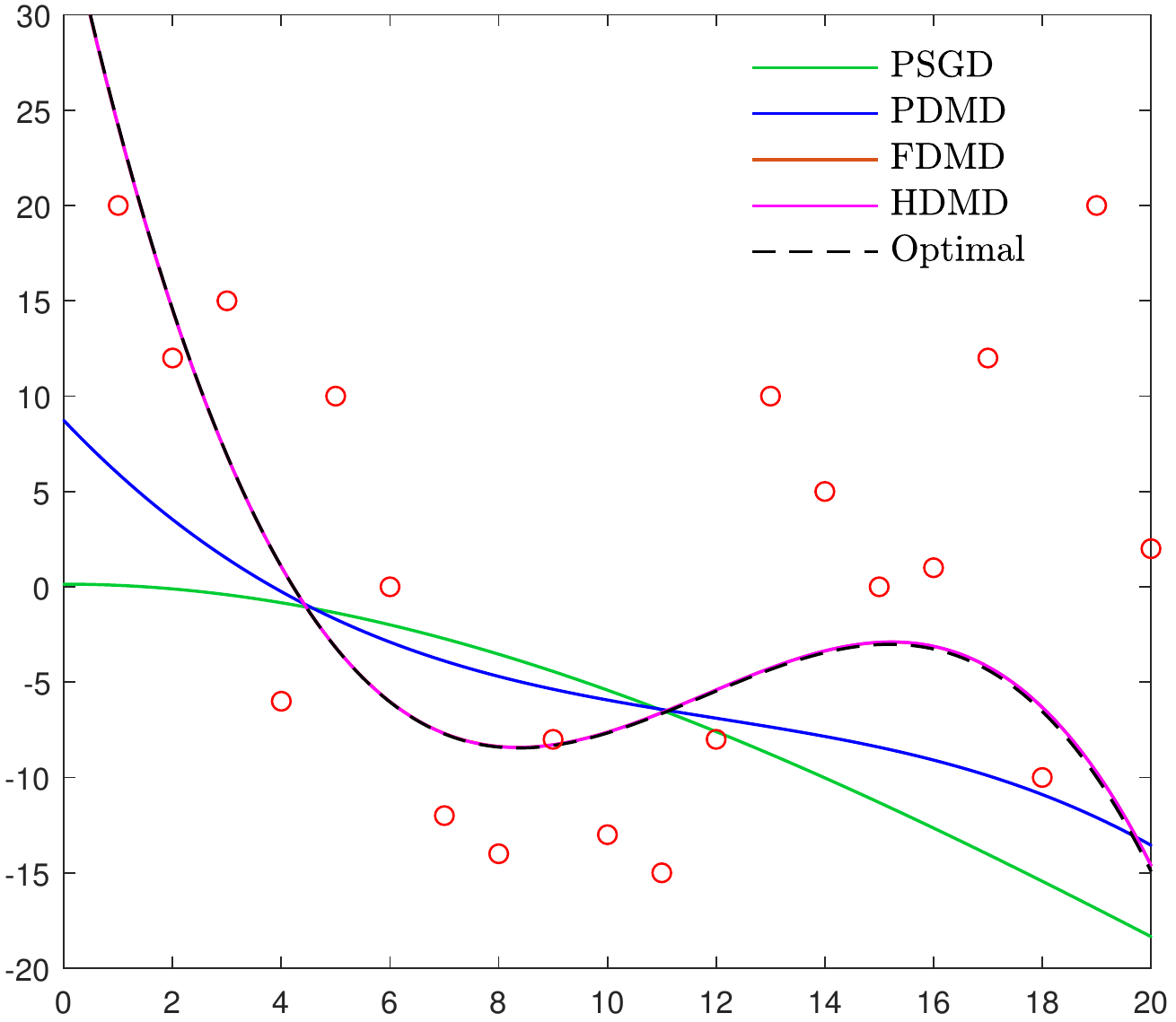}}
	\caption{Final third-order fit}   
	\label{fig:Polynomial_Regression_Curve_Fit_Third_Order_Epsilon_0_1}
\end{minipage}
\end{figure}


\end{example}

\vspace{-0.7cm}
\section{Conclusions}
In this paper, we proposed two continuous-time classes of mirror dynamics for monotone concave games. 
 We showed that they are guaranteed to converge to a perturbed Nash equilibrium, which  tends to a Nash equilibrium of the game 
as the regularization goes to zero. 
One the two classes allows convergence in hypo-monotone games. 
We provided several examples from both classes. 
As future work, we will consider incomplete information, where players observe only a portion of the partial-gradient or a noisy version of it. 

	\vspace{-0.35cm}
\appendix

	\vspace{-0.35cm}

\begin{proof}\emph{ of  \autoref{prop:steep_implies_differentiability_in_interior}}	
\begin{itemize}
		\item[(i)] Suppose that $\dom(\vartheta^p)\! \!\subset\! \!\mathbb{R}^{n_p}$ has dimension $n_k \!\!< \!\!n_p$, then there exists a one-to-one affine transformation $T(x) \!\!=\!\! Mx \!\!+\!\! b$, $M: \mathbb{R}^{n_p} \!\to \! \mathbb{R}^{n_p}$ linear, of $\mathbb{R}^{n_p}$ onto itself which maps $\aff(\dom(\vartheta^p))$ onto the subspace $\mathcal{L} \!=\! \{x^p \!=\! (x^p_1, \ldots, x^p_{n_k}, x^p_{n_k+1}, \ldots x^p_{n_p})|x^p_{n_k+1} \!=\! 0, \ldots, x^p_{n_p} \!=\! 0\}$ \cite[p. 45]{Rockafellar}. Then $T(\dom(\vartheta^p))\!\! \subset \!\!T(\aff(\dom(\vartheta^p))) \!\!=\!\! \mathcal{L}$. Since the subspace $\mathcal{L}$ is homeomorphic to $\mathbb{R}^{n_k}$, therefore $\interior( T(\dom(\vartheta^p)))$ is non-empty when regarded as a subset of $\mathbb{R}^{n_k}$. Then the result follows by applying Theorem 26.1 of \cite[p. 251]{Rockafellar} to $\widetilde \vartheta^p(x^p_1, \ldots, x^p_{n_k}) \!\!=\!\! \vartheta^p(x^p_1, \ldots, x^p_{n_k}, 0, \ldots, 0)$, the restriction of $\vartheta^p$ to $\mathcal{L}$. Otherwise if the dimension of $\dom(\vartheta^p)$ is $n_p$, then the interior  coincides with the relative interior, the result again follows from Theorem 26.1 of \cite[p. 251]{Rockafellar}.
		\item[(ii)]  $\rboundary(\dom(\vartheta^p)\! \!\subset\!\! \dom(\partial \vartheta^p)$ follows from the definition. 
		By \cite[p. 227]{Rockafellar}, $\rinterior(\dom(\vartheta^p))\! \!\subseteq \! \!\dom(\partial \vartheta^p) \!\!\subseteq\! \dom(\vartheta^p)$. Since $\rinterior(\dom(\vartheta^p)) \!\!\cap \!\!\rboundary(\dom(\vartheta^p)) \!\! =\!\! \varnothing$, $\rinterior(\dom(\vartheta^p)) \!\!\cup\!\! \rboundary(\dom(\vartheta^p)) \!\!=\!\! \closure(\dom(\vartheta^p))$, then $\closure(\dom(\vartheta^p)) \! \!\subseteq \! \!\dom(\partial \vartheta^p)$, hence 
		$\dom(\vartheta^p) \!\! \subseteq \! \!\dom(\partial \vartheta^p)$, which shows the reverse. 
	\end{itemize}\vspace{-0.4cm}\end{proof}
	\vspace{-0.5cm}
\begin{proof}\emph{of \autoref{prop:primal_sc}}
	\begin{enumerate}
		\item[(i)] ${\psi^p}^\star$ is closed, proper, convex follows from \cite[Theorem 4.4, p. 87]{Beck17}, \cite[Theorem 4.5, p. 88]{Beck17}. Since $\psi^p$ is closed, proper, $\epsilon\rho$-strongly convex, therefore it is supercoercive \cite[Prop 3.10.8, p. 169]{Niculescu}, hence 
		${\psi^p}^\star$ is finite for all $z^p \!\! \in \!\! \mathbb{R}^{n_p}$. 
		\item[(ii)] The continuous differentiablility of ${\psi^p}^\star$ follows from \cite[Theorem 6.2.4(a), p. 264]{Niculescu}. Since $\psi^p$ is strongly convex and supercoercive, therefore the maximizer of \eqref{eqn:mirror_map_argmax_char} exists, is unique and equals to $\nabla {\psi^p}^\star$ for all $z^p \in \mathbb{R}^{n_p}$,  hence $C^p \!=\! \nabla {\psi^p}^\star$. 
		\item[(iii)] Since $\psi^p$ is closed, proper, $\epsilon\rho$-strongly convex, the Lipschitzness of $C^p$ follows from \cite[Theorem 5.26, p. 123]{Beck17}. 
		\item[(iv)] Since $C^p \!=\! \nabla {\psi^p}^\star\!$ and $C^p$ is $(\epsilon\rho)^{-1}$-Lipschitz by (iii), $\epsilon\rho$-cocoercivity follows from the Baillon-Haddad theorem, see \cite[Corollary 18.17, p. 323]{Bauschke2nd}, and monotonicity directly follows. 
		\item[(v)] Since ${\psi^p}^\star$ is proper, closed and convex,  by \cite[Theorem 4.20, p. 104]{Beck17}, $\partial {\psi^p}^\star$ is an inverse of $\partial \psi^p$ and vice-versa, i.e., $x^p \!\! \in \!\! \partial {\psi^p}^\star(z^p)$ if and only if $z^p \!\! \in \!\! \partial \psi^p(x^p)$.  Since $\partial {\psi^p}^\star(z^p) \!\!= \!\!\{\nabla {\psi^p}^\star \!(z^p)\}, \!\forall z^p \in \mathbb{R}^{n_p}$, then $x^p \!\!=\!\! \nabla {\psi^p}^\star \!(z^p)$. Let $x^p$ be such that $\partial \psi^p(x^p)\!\! \neq \!\!\varnothing$, then the set of all such $x^p$ is  the domain of $\partial \psi^p$, hence $\range(\partial {\psi^p}^\star) \!\!=\!\! \range(\nabla {\psi^p}^\star) \!\!=\!\! \range(C^p) \subseteq \dom(\partial \psi^p)$. By a similar argument,  $ \dom(\partial \psi^p)\! \!\subseteq \! \!\range(C^p)$, so $C^p$ is surjective from $\mathbb{R}^{n_p}$ onto $\!\dom(\partial \! \psi^p)$. 
		If $\psi^p$ is steep, by \autoref{prop:steep_implies_differentiability_in_interior}, $\dom(\partial \psi^p) \!\! =\!\! \rinterior(\dom(\psi^p))$, hence $\range(C^p) \!\!=\! \!\rinterior(\dom(\psi^p))$. Otherwise, if $\psi^p$ is non-steep, by \autoref{prop:steep_implies_differentiability_in_interior}, $\dom(\partial \psi^p) \!\!=\!\! \dom(\psi^p)$, hence $\range(C^p) \!\!=\!\! \dom(\psi^p) $.

		\item[(vi)] From (v), since $x^p \!\! \in  \!\!\partial {\psi^p}^\star \!(z^p)$ if and only if $z^p \!\!\in \!\!\partial \psi^p(x^p)$, and $x^p \!\! =\!\! \nabla {\psi^p}^\star \!(z^p) \! =\! C^p(z^p)$, therefore $z^p \!\!\in 
		\! \!\partial \psi^p(C^p(z^p))$. By (v),  $C^p(z^p) \!\!\in \!\!\dom(\partial \psi^p)$, hence $\partial \psi^p$ is a singleton,  $\partial \psi^p\! =\! \{\nabla \psi^p\}$, 
		so that $z^p\! \!=\!\! \nabla \psi^p(C^p(z^p)\!)$, i.e. $\nabla \psi^p$ is a left-inverse of $C^p$.	\end{enumerate}\vspace{-0.4cm}\end{proof}
	\vspace{-0.5cm}
\begin{proof}\emph{of \autoref{prop:mirrormap_sc}}
	Consider the argmax characterization \eqref{eqn:mirror_map_argmax_char} for $\!\psi^p(y^p) \! \!=\! \! \epsilon \vartheta^p(y^p)$. Thus, 
	$C^p(z^p) 
	\!\!=\! \! \underset{y^p \in \Omega^p}{\text{argmin}}\left[\psi^p(y^p) \!-\!  {z^p}^\top y^p \right] 
	 \!\!= \!\! \underset{y^p \in \Omega^p}{\text{argmin}}\left[\psi^p(y^p) \! -\! \psi^p(\nabla {\psi^p}^\star(z^p))\! -\!  {z^p}^\top (y^p \!-\! \nabla {\psi^p}^\star(z^p))\right] 
	 $,  	
	by inserting terms independent of $y^p\!$. Under  \autoref{assump:primal}(i), by \autoref{prop:primal_sc}(v), $\nabla \psi^p$ is a left-inverse of $\nabla {\psi^p}^\star\!$ over $\range(C^p)\! $ $\!=\!\!\rinterior(\dom(\psi^p\!))$ if $\psi^p$  is steep, and over  
	$ \!\dom(\psi^p)\! \!=\!\! \Omega^p\!$ if $\psi^p\!$  is non-steep. Alternatively, under \autoref{assump:primal}(ii), by \autoref{prop:primal_Legendre},  $\nabla \psi^p\!$ is inverse of $\nabla {\psi^p}^\star\!$ over $\range(C^p) \!\!=\!\! \interior(\dom({\psi^p})\!) \!\!$ $\!=\!\! \interior(\Omega^p)$. Therefore, $z^p \! \!=\!\! \nabla \psi^p(\nabla {\psi^p}^\star\!(z^p))$, which used in the last term 
yields the Bregman divergence of $\psi^p$. 
\vspace{-0.4cm}\end{proof}

\vspace{-0.2cm}

\end{document}